\documentclass{article}

\usepackage[english]{babel}

\usepackage[    
    style=numeric-comp,    % citazioni numeriche [1], [2], ...
    maxnames=6,       % se più di 3 autori, usa et al.
    minnames=5,
    sortcites=true,      % ordina i numeri nelle citazioni multiple
    doi=false,        % niente DOI
    isbn=false,       % niente ISBN
    url=false,        % niente URL
    eprint=true,      % arXiv/eprint
    giveninits=true   % abbreviate first names
    ]
{biblatex}
\usepackage[letterpaper,top=2cm,bottom=2cm,left=3cm,right=3cm,marginparwidth=1.75cm]{geometry}
\usepackage{authblk}
\usepackage{amsmath}
\usepackage{graphicx}
\usepackage[colorlinks=true, allcolors=blue]{hyperref}
\usepackage{subfigure}
\usepackage{todonotes}
\usepackage{amssymb}

\newcommand{\RR}{\mathbb{R}}
\usepackage{authblk}

\date{}

\title{High-Order Asymptotic-Preserving IMEX schemes for an ES-BGK model for Gas Mixtures}

\author[1,2]{Domenico Caparello\thanks{Domenico.caparello@univ-cotedazur.fr, domenico.caparello@unife.it}}
\author[2,3]{Lorenzo Pareschi\thanks{L.Pareschi@hw.ac.uk}}
\author[1]{Thomas Rey\thanks{Thomas.rey@univ-cotedazur.fr}}

\affil[1]{Université Côte d’Azur, CNRS, LJAD, Parc Valrose, F-06108 Nice, France}
\affil[2]{Department of Mathematics and Computer Science University of Ferrara}
\affil[3]{Heriot-Watt University, Edinburgh, UK}

\begin{document}

\maketitle

\begin{abstract}
In this work we construct a high-order Asymptotic-Preserving (AP) Implicit-Explicit (IMEX) scheme for the ES-BGK model for gas mixtures introduced in \cite{bru-2015}. The time discretization is based on the IMEX strategy proposed in \cite{fil-jin-2011} for the single-species BGK model and is here extended to the multi-species ES-BGK setting. The resulting method is fully explicit, uniformly stable with respect to the Knudsen number and, in the fluid regime, it reduces to a consistent and high-order accurate solver for the limiting macroscopic equations of the mixture. The IMEX structure removes the stiffness associated with the relaxation term so that the time step is constrained only by a hyperbolic CFL condition. 
The full solver couples a high-order space and velocity discretization that includes third-order time integration, a CWENO3 finite-volume reconstruction in space, exact conservation of macroscopic moments in the discrete velocity space, and a multithreaded implementation. The proposed approach can handle an arbitrary number of species. Its accuracy and robustness are demonstrated on a set of multidimensional kinetic tests for gas mixtures, where the AP property and the correct asymptotics are numerically verified across different regimes.
\end{abstract}

\section{Introduction}
Mixtures of rarefied gas can be found in a wide variety of systems, such as high-altitude gas, nuclear engineering applications, evaporation-condensation processes \cite{jin-qin-2013}, or in the context of plasma physics, by considering, for example, the evolution of electrons and ions as different species \cite{bot-cow-sch-2024}. Here we describe some examples (coming from \cite{jin-qin-2013}) in which macroscopic models (like Euler or Navier-Stokes equations) fail to correctly describe the fluid, and thus, kinetic models of the Boltzmann types are needed to describe these systems. In the space shuttle reentry problem, at high altitudes, free streaming effects are significant and thus macroscopic models are not adequate in describing the system evolution. The shock profiles computed by using the Navier-Stokes equations are not accurate for hypersonic flows, while the Boltzmann equation is able to reproduce the correct behaviour.  Thus, given the complexity of gas mixtures, a realistic mathematical depiction  of such systems must involve multi-species kinetic models.

The analytical formulation of the Boltzmann equation for gas mixtures has been examined in several works, like \cite{mar-2016, bri-dau-2016}. In \cite{rey-ten-2024} the authors extend the gas mixtures Boltzmann equation framework to granular gases and study its well posedness in the homogeneous case as well as its large time behavior.

Such kind of microscopic kinetic models (like the Boltzmann equation) often contains small parameters, like for example the mean free path of particles, and an asymptotic expansion of the kinetic models on these parameters typically yield macroscopic equations (like the Euler equations), which describe the average behaviour of the fluid. The asymptotic resolution of these scales using kinetic models gives rise to significant numerical challenges due to the high computational cost required, thus in case of validity of the macroscopic equations, it is much less expensive to use them instead of integrating the microscopic kinetic model. However, we have already seen that in some cases such macroscopic equations can break down in part of the domain, and thus kinetic microscopic models are needed in such cases. A possible approach for dealing with these scenarios is by using multi-physics dynamic domain decomposition strategies, where the scheme employs macroscopic or kinetic equations according to the non-equilibrium behaviour of the fluid, see for example \cite{fil-rey-2015, tiw-1998, cap-par-rey-2025}. Another approach for such multiscale problems, which we also adopt in our scheme, is to employ an asymptotic-preserving (AP) strategy. The core idea behind these schemes is to ensure that the macroscopic models are correctly recovered at the discrete level in the asymptotic limit of the underlying microscopic model. Furthermore, for a scheme to be classified as asymptotic-preserving it must fulfil a second condition, namely the implicit collision term must be implemented efficiently, typically allowing the time step to be much larger than the \emph{small parameter} (typically the Knudsen number, i.e. the ratio between the mean free path of particles and the physical length scale of the system) of the problem. This brief introduction to AP schemes that we have just presented is based on \cite{shi-2010}.
Several works focused on Asymptotic-Preserving schemes for kinetic equation can be found in literature, as \cite{diaz-chen-2015, fil-jin-2011, fil-jin-2010, jin-qin-2013, hu-wan-2015, bos-cho-rus-yun-2021}.\\

Gas mixtures models describe the evolution of the system by evolving the distribution of particles $f=f(t, x, v)$ of each species, which is a non-negative function representing the probability to find a particle in a volume $dx$ around $x$, and in an interval $dv$ around $v$, during a small time $dt$. The distribution function is normalized such that $\int \int f(t,x,v) \  dx \  dv=N$ where $N$ is the total number of particles. Mixture of $P$ species are each described by a particle distribution function $f_p(t,x,v)$ for $p \in \{1,\ldots, P\}$ and normalized such that $\int \int f_p(t,x,v) \  dx \  dv=N_p$, where $N_p$ is the total number of particles of the $p$-th species.

Given $t\geq 0$, $x\in \Omega \subset \RR^{D_X}$ and $v\in \RR^{D_V}$, each function $f_p$ evolves according the multi-species Boltzmann equation
\begin{equation}
    \begin{cases}
    \displaystyle \partial_t f_p + v \cdot \nabla_xf_p = \frac{1}{\varepsilon}\mathcal{Q}_p[\mathbf{f}],\\
    \displaystyle \mathcal{Q}_p[\mathbf{f}]=\sum_{q=1}^P \mathcal{Q}_{p,q}[f_p, f_q],\\
    f_p(0, x, v)=f_{p,0}(x,v),   
    \end{cases}
\end{equation}
where $\mathcal{Q}_P[\mathbf{f}]$ is the collision operator which describes the changes in momentum due to collisions between particles and $\varepsilon$ is the Knudsen number, which is the ratio between the mean free path of particles and the physical length scale of the system: thus, if $\varepsilon<<1$ then the fluid can be described using the so-called \emph{hydrodynamic limit} (for example, using the compressible Euler equations, Navier-Stokes or Burnett, etc.).

The collision operator is given by \cite{bai-rey-2021}
\begin{equation}
    \mathcal{Q}_{p,q}[f_p, f_q] = \int \int_{\RR^{D_V}\times \mathbb{S}^{D_V - 1}}B_{p,q}\left(|v-v_*|, \cos{\theta}\right) \left[f_p(v')f_q(v_*')-f_p(v)f_q(v_*)\right]\  dv_*\  d\sigma,
\end{equation}
where the inter-species collisions can be parametrized as %velocities of the colliding pairs $\left(v, v_*\right)$ and $\left(v', v_*'\right)$ can be parametrized as
\begin{equation}
    \begin{cases}
        v'=\displaystyle\frac{1}{m_p + m_q}\left(m_p v+m_q v_*+m_q|v-v_*|\sigma\right),\\
        v'_*=\displaystyle\frac{1}{m_p + m_q}\left(m_p v+m_q v_*-m_p|v-v_*|\sigma\right),
    \end{cases}
\end{equation}
where $m_p$ is the mass of each particle of the $p$-th species.\\
The quantity $B_{p,q}$ is the inter-species collision kernel: a non-negative function describing the physics of the collisions between particles. More details on the collision kernel for the single species case (and for the Boltzmann equation in general) can be found in \cite{dim-par-2014}.\\
The collision operator must preserve the mass of each species, the total momentum and the total kinetic energy of the system. In formulas, the following quantities must be conserved
\begin{equation}
    \begin{cases}
        \displaystyle \int f_p \  dv, \quad \quad \quad \quad \forall p \in \{1, 2, ..., P\},\\
        \displaystyle \sum_{p=1}^P \int \int m_p v f_p(t,x,v) \  dx \  dv,\\
        \displaystyle \sum_{p=1}^P \int \int \frac{m_p}{2} |v|^2 f_p(t,x,v) \  dx \  dv.\\
    \end{cases}
\end{equation}
Every species has its own global equilibrium distribution, that are all Maxwellians with common average velocity $u$ and temperature $T$
\begin{equation}\label{eq7}
    \mathcal{M}_p^{n_p, u, T}(t,x,v)=
    n_p(t,x)\left(\frac{m_p}{2\pi T}\right)^{\frac{D_V}{2}}
    \exp{\left(-\frac{m_p (u-v)^2}{2T}\right)},
\end{equation}
where the numerical density $n_p$, the average velocity $v_p$ and the temperature $T_p$ of the $p$-th species are given by
\begin{equation}
    n_p=\int f_p \  dv, \quad \quad \quad
    \rho_p v_p = \int m_p  v f_p \  dv, \quad \quad \quad
    \frac{D_V}{2}n_p T_p = \int \frac{m_p}{2}(v_p-v)^2 f_p \  dv,
\end{equation}
with $\rho_p=m_p n_p$ and  $\mathcal{P}_p=n_p T_p$; these last two quantities are, respectively, the density and the pressure of the $p$-th species, and $D_V$ is the dimensions of the velocity space. The analogous total moments for the mixture are given by $n=\sum n_p$, $\rho=\sum \rho_p$, $\rho u = \sum \rho_p v_p$, $\mathcal{P}=\sum\mathcal{P}_p$ and
\begin{equation}
    \frac{D_V}{2}nT = \sum_{p=1}^P \int \frac{m_p}{2}(u-v)^2 f_p \  dv.
\end{equation}
It is possible to define also a macroscopic total energy of the mixture as follows
\begin{equation}
    E= \sum_p\int \frac{m_p}{2} |v|^2f_p dv=\frac{1}{2}\rho |u|^2+\frac{D_V}{2}nT.
\end{equation}
In literature several possible models for the multi-species Boltzmann equation have been proposed. A possible approach consists in a generalization of the single species BGK model \cite{bha-gro-kro-1954}: the main idea is that the collision operator drives the distribution function towards the equilibrium Maxwellian through a simple relaxation mechanism. As discussed in \cite{gro-mon-spi-2011}, the extension of this BGK strategy to gas mixtures is not trivial, since it can give rise to breakdown of positivity for density or temperature fields, and/or of the indifferentiability principle.
Another possible approach is to replace the collision operator with a sum of relaxation operators (like \cite{haa-hau-mur-2017, %kli-pir-pup-2017, %mcc-1973,
bob-bis-gro-spi-pot-2018}) of the BGK kind, each describing the interactions with any other species (itself or different).\\
In \cite{bis-bos-dim-gro-mar-2022} the authors propose a binary mixture model where the collisions between particles of the same species are treated using the full Boltzmann operator, while the BGK relaxation terms are used to describe the interactions between particles of different species. They also validate their model with some interesting numerical simulations.
In the recent paper \cite{zha-zha-son-guo-2025}, the authors propose an implicit discrete unified gas kinetic scheme (DUGKS) for multiscale steady flows of binary gas mixtures.\\
Of particular interest is \cite{bos-cho-gro-rus-2021}, where the authors present a comparison of different models for gas mixtures proposed by \cite{and-aok-per-2002, bis-gro-spi-2010, bob-bis-gro-spi-pot-2018}. In \cite{jin-qin-2013} the authors extends to the multi-species Boltzmann equation the BGK penalization technique of the full Boltzmann operator introduced by \cite{fil-jin-2010}: this technique allows to rewrite the Boltzmann equation as a sum of a stiff and a non-stiff term.\\
In \cite{gro-mon-spi-2011} the authors propose an ellipsoidal BGK model for a binary mixture in the kinetic theory framework. Different spectral methods for the full multi-species Boltzmann operator have been proposed, for example, by \cite{wu-zha-zha-ree-2015, sha-ali-jin-2019}.\\
In this context, the derivation of the corresponding hydrodynamic equations from kinetic models \cite{lev-1996} is of particular interest, for example for inert gas mixtures, reactive mixtures or multi-phase flows \cite{bis-gro-mar-2021, pup-rey-ten-2025, bis-gro-spi-2010}.\\

A big challenge in numerical schemes for the multi-species Boltzmann equation (as well as in plasma contexts) arises from the non-unity mass ratios among species. This difficulty becomes more problematic as the mass ratio increases, since larger velocity domains are required, which in turn necessitates more discretization cells and consequently leads to higher computational costs \cite{bai-rey-2021}. Our simulations will show that our numerical scheme is able to correctly integrate unequal-mass configurations. In particular we perform some tests with a mass ratio equal to 100 in the spatially homogeneous settings, and up to 20 in the spatially non-homogeneous setting.\\
In this work we are interested in developing an AP-IMEX scheme for the  collision operator presented by Brull in \cite{bru-2015}
\begin{equation}
    \mathcal{Q}_p = \lambda \left(G_p - f_p\right),
\end{equation}
where $\lambda$ is a free parameter that can be used to fit the transport coefficients, and
\begin{equation}\label{eq2000}
    \displaystyle G_p = \frac{n_p}{\sqrt{\det{(2\pi\tilde{\tau}_p)}}}
    \exp{\left(-\frac{1}{2}(v-u)^\intercal\tilde{\tau}_p^{-1}(v-u)\right)},
\end{equation}
and the corrected tensors
\begin{equation}\label{eq3}
    \begin{cases}
    \displaystyle \tilde{\tau}_p =\frac{1}{n_p} \int (v-u)\otimes(v-u)\  G_p \  dv =  \tau \frac{\rho}{n} \frac{1}{m_p} ,\\
    \displaystyle \tau = \nu \Theta + (1-\nu)T \frac{n}{\rho} \text{I},    
    \end{cases}    
\end{equation}
where I is the identity matrix, the parameter $\nu$ is used to modify the value of the Prandtl number for $\nu \in [-1/2, 1)$ and the following quantity has been introduced
\begin{equation}\label{eq6000}
    \Theta = \frac{1}{\rho}\sum_{p=1}^P \int m_p (v-u)\otimes(v-u) \  f_p \  dv.
\end{equation}
This operator \cite{bru-2015} is an extension of the Ellipsoidal Statistical Model for the single species Boltzmann equation \cite{%hol-low-1965, 
cer-1988, low-hol-1966}, and its derivation is based on the resolution of a minimization problem, given some relaxation coefficients associated to some moments \cite{bru-pav-sch-2012, bru-sch-2008, bru-sch-2009}. In particular, in \cite{bru-sch-2008} the authors propose a new approach of the ES-BGK model based on an entropy minimization principle \cite{jun-2000, lev-1996, sch-2004}: different relaxation rates related to some moments are defined, and they describe the way in which these moments vanish for the trend to equilibrium. This approach has been applied to polyatomic gases in \cite{bru-sch-2009}.\\

The AP property of our scheme, which is an extension of the work of \cite{fil-jin-2011}, is based on an Implicit-Explicit formulation of the problem. This particular structure, allows to solve explicitly both the transport term and the collision term, despite having discretized this last term implicitly. Thanks to the implicit structure, the choice of the timestep is not restricted anymore by the choice of the \emph{small parameter} of the problem (the Knudsen number in our case), but only by the CFL transport condition. This particular formulation allows for a significant reduction in computational costs (as will be discussed in more detail at the end of the paper), since in stiff regimes the Knudsen number may become very small. In such cases, a fully explicit solver would be subject to a much more restrictive stability condition than the standard transport CFL constraint, as the timestep cannot be larger than the Knudsen number as order of magnitude. Consequently, a fully explicit scheme would result in a prohibitive large computational cost.\\

To construct a high-order solver for the two-dimensional ES-BGK system for gas mixtures in both physical and velocity space, we combine several advanced numerical ingredients into a coherent and fully Asymptotic-Preserving framework. \begin{itemize}
    \item a third-order IMEX Runge--Kutta time integrator for stiff kinetic equations \cite{par-rus-2001, bos-par-rus-2024};
    \item a high-resolution CWENO3 finite-volume reconstruction that provides third-order spatial accuracy \cite{Lev-Pup-Rus-1999, pup-sem-vis-2023};
    \item a conservative velocity discretization that preserves the macroscopic moments of each species, obtained by extending the weighted $L^2$ minimization technique of \cite{bos-cho-rus-2022} to the multi-species setting;
    \item a multithreaded CPU implementation that ensures computational efficiency in high-dimensional phase space;
    \item the Asymptotic-Preserving IMEX formulation specifically developed in this work, which generalizes the construction of \cite{fil-jin-2011} to the multi-species ES-BGK model of \cite{bru-2015}.
\end{itemize}
To the best of our knowledge, this is the first work that combines a high-order AP IMEX time discretization with a multi-species ES-BGK collision operator in multiple space dimensions, together with a high-order spatial and velocity discretization and a parallel implementation capable of handling an arbitrary number of species. 

The rest of the manuscript is organized as follows. In Section \ref{secSchemes} we introduce the AP IMEX formulation based on \cite{fil-jin-2011} for the ES-BGK model for gas mixtures proposed by \cite{bru-2015}, we also prove the Asympotic-Preserving property of the method, we discuss how to arbitrarily increase the temporal order, and then we present the techniques that we implemented for the conservation of macroscopic moments. In Section \ref{secNumTest} we show some numerical tests used to validate the scheme, and in the last Section \ref{conclusions} we summarize the key aspects of this work.

\section{An asymptotic-Preserving scheme for the ES-BGK equation for gas mixtures}\label{secSchemes}
As summarized in \cite{fil-jin-2011, jin-1999, shi-2010} a scheme for the Boltzmann equation can be defined AP (Asymptotic-Preserving) if
\begin{itemize}
    \item it preserves the discrete analogy of the Chapman-Enskog expansion: when holding the mesh size and time step fixed and letting the Knudsen number goes to zero, then one recovers a suitable scheme for the limiting Euler equations;
    \item implicit collision is treated explicitly, or at least, implemented in a manner that is computationally more efficient than the use of Newton-type solvers for nonlinear algebraic systems.
\end{itemize}
\subsection{A first-order IMEX strategy}\label{subsectionIMEX}
The following AP-IMEX scheme for this multi-species ES-BGK model is strongly inspired by the one proposed by Filbet and Jin in \cite{fil-jin-2011}, and we will present it in a three-dimensional velocity space.\\
For each species, a first order IMEX scheme for the Boltzmann equation with this multi-species ES-BGK operator will read as
\begin{equation}\label{eq1002}
    \begin{cases}
        \displaystyle \frac{f_p^{n+1}-f_p^n}{\Delta t} + v \cdot \nabla_x f_p^n = \frac{\lambda^{n+1}}{\varepsilon}\left(G_p^{n+1} - f_p^{n+1}\right),\\
        \displaystyle f_p^0(x,v) = f_p(0,x,v)
    \end{cases}
\end{equation}
The previous equation can be rearranged as
\begin{equation}\label{eq5}
    f_p^{n+1}=\frac{\varepsilon}{\varepsilon + \lambda^{n+1}\Delta t}\left[f_p^n-\Delta t \  v\cdot \nabla_x f_p^n\right] + \frac{\lambda^{n+1}\Delta t}{\varepsilon + \lambda^{n+1}\Delta t} G_p^{n+1}.
\end{equation}

Integrating each of these last $P$ equations with respect to velocity, and using the conservative properties of the $\mathcal{Q}_p$ then it is possible to compute the numerical density of each species at the next timestep from the transport term. Analogously, multiplying each of these last $P$ equations by $m_p \ v$ and by $m_p \ |v|^2$, integrating with respect to velocity, summing up all the equations with respect to the species and using the conservative properties of $\mathcal{Q}_p$ then it is possible to compute the corresponding total macroscopic quantities of the fluid (average velocity and average temperature). Doing that we obtain the following equations for $n_p^{n+1}$, $u^{n+1}$ and $E^{n+1}$ (and thus $T^{n+1}$) as follows
\begin{equation}\label{eq4}
    \begin{cases}
        \displaystyle n_p^{n+1} = \int \left(f_p^n - \Delta t \  v\cdot \nabla_x f_p^n\right)\  dv,\\
        \displaystyle \rho^{n+1}=\sum_{p=1}^P \int m_p \left(f_p^n - \Delta t \  v\cdot \nabla_x f_p^n\right) \  dv=\sum_{p=1}^P m_p n_p^{n+1},\\
        \displaystyle \rho^{n+1} u^{n+1} = \sum_{p=1}^P \int m_p v \left(f_p^n - \Delta t \  v\cdot \nabla_x f_p^n\right)\  dv,\\
        %\displaystyle \frac{D_V}{2}n^{n+1}T^{n+1} = \sum_{p=1}^P \int \frac{m_p}{2}(u^{n+1}-v)^2 \left(f_p^n - \Delta t v\cdot \nabla_x f_p^n\right) dv.
        \displaystyle E^{n+1} = \sum_{p=1}^P \int \frac{m_p}{2}|v|^2 \left(f_p^n - \Delta t\  v\cdot \nabla_x f_p^n\right) dv,
    \end{cases}
\end{equation}
and using the relation
$E=\sum_p \int  \frac{m_p}{2} |v|^2f_p \  dv=\frac{1}{2}\rho |u|^2+\frac{3}{2}nT$ then it is possible to compute $T^{n+1}$.\\

Unfortunately, this is not sufficient to compute the operators $G_p^{n+1}$.\\
So let's introduce the following quantity
\begin{equation}\label{eq2}
    \displaystyle \Sigma^{n+1} := \sum_{p=1}^P \int v\otimes v \  m_p f_p^{n+1} \  dv=\rho^{n+1}\left[\Theta^{n+1} + u^{n+1}\otimes u^{n+1}\right],
\end{equation}
where we used equation \eqref{eq6000} in this last equality.\\
Let's note that
\begin{equation}\label{eq6}
    \sum_{p=1}^P \int m_p (v\otimes v) G_p  \  dv = \rho \left( \tau + u\otimes u\right).
\end{equation}
So multiplying equation \eqref{eq5} by $m_p \  v\otimes v$, integrating with respect to velocity, summing up with respect to all species, and using equations \eqref{eq2} and \eqref{eq6} then one gets
\begin{equation}
    \begin{aligned}
        \displaystyle \Sigma^{n+1} = \frac{\varepsilon}{\varepsilon + \lambda^{n+1}\Delta t}
        \left(\Sigma^n - \Delta t\sum_{p=1}^P \int  m_p v\otimes v \ v\cdot\nabla_x f_p^n \  dv\right)\\
        +\frac{\lambda^{n+1}\Delta t}{\varepsilon+\lambda^{n+1}\Delta t }\rho^{n+1} \left(\tau^{n+1}+u^{n+1}\otimes u^{n+1}\right).
    \end{aligned}    
\end{equation}
Using equation \eqref{eq3}
\begin{equation}
    \begin{aligned}
        \displaystyle \Sigma^{n+1} = \frac{\varepsilon}{\varepsilon + \lambda^{n+1}\Delta t}
        \left(\Sigma^n - \Delta t\sum_{p=1}^P \int  m_p v\otimes v \  v\cdot\nabla_x f_p^n \  dv\right)\\
        +\frac{\lambda^{n+1}\Delta t}{\varepsilon+\lambda^{n+1}\Delta t }\rho^{n+1} \left(\nu \Theta^{n+1} + (1-\nu) T^{n+1}\frac{n^{n+1}}{\rho^{n+1}}\text{I}+u^{n+1}\otimes u^{n+1}\right).
    \end{aligned}    
\end{equation}
The quantity $\Theta^{n+1}$ can be replaced by $\Sigma^{n+1}/\rho^{n+1} - u^{n+1}\otimes u^{n+1}$ (see equation \eqref{eq2}), and by factoring out the therm in $\Sigma^{n+1}$ on the left-hand side the following equation is obtained
\begin{equation}\label{eq1}
    \begin{aligned}
        \displaystyle \Sigma^{n+1} = \frac{\varepsilon}{\varepsilon + \lambda^{n+1}\Delta t(1-\nu)}
        \left(\Sigma^n - \Delta t\sum_{p=1}^P \int  m_p v\otimes v \  v\cdot\nabla_x f_p^n dv\right)\\
        +\frac{\lambda^{n+1}\Delta t (1-\nu)}{\varepsilon+\lambda^{n+1}\Delta t (1-\nu)}\left[T^{n+1}n^{n+1}\text{I}  + \rho^{n+1} u^{n+1}\otimes u^{n+1}\right].
    \end{aligned}    
\end{equation}
The quantity in round brackets on the right-hand side can be directly computed since the distribution functions $f_p^n$ are known, while the macroscopic moments in square brackets on the right-hand side can be computed explicitly from $f_p^n$ using equation \eqref{eq4}. By virtue of equation \eqref{eq1} it is possible to compute $\Sigma^{n+1}$, and then the quantity $\Theta^{n+1}$ by using equation \eqref{eq2}. After that $\tau^{n+1}$ and $\tilde{\tau}_p^{n+1}$ are obtained thanks to equation \eqref{eq3} and in the end $G_p^{n+1}$ (equation \eqref{eq2000}).\\
At this point it is possible to update each distribution function by computing its value at the next timestep by using equation \eqref{eq5}. It is important to remark that this IMEX strategy (based on \cite{fil-jin-2011}) allows to solve \emph{explicitly} the equation \eqref{eq1002} which is nonlinearly implicit. In the following section we will discuss per AP property of this scheme. 
\subsection{Asymptotic-Preserving property}\label{APsection}
Now we prove the Asymptotic-Preserving property of the scheme similarly to \cite{fil-jin-2011}.\\
\textbf{Proposition 1}. Consider the solution $f_p^{n+1}$ computed as indicated in the last section \ref{subsectionIMEX} then,
\begin{enumerate}
    \item For all $\varepsilon\to 0$ and $\Delta t>0$, the distribution functions $f_p^{n+1}$ satisfy $$ 0\leq f_p^{n+1}(x,v)\leq \max{(||f^n_p||_\infty, ||G_p^{n+1}||_\infty)}.$$
    \item For all $\Delta t >0$ and $f_p^0$, the distribution functions $f_p^n$ converge to the equilibrium Maxwellians $\mathcal{M}_p^{n_p^n, u^n, T^n}$, that is $$\lim_{\varepsilon \to 0}f_p^n=\mathcal{M}_p^{n_p^n, u^n, T^n},$$ and the scheme gives a first order approximation in time of the compressible Euler system for gas mixtures.
    \item Assuming that $||f_p(t^n) - \mathcal{M}_p(t^n)|| = \mathcal{O}(\varepsilon)$, for $n\geq 2$ and
    \begin{equation}\label{eqAssumption}
        \left|\left|\frac{U^{n+1} - U^n}{\Delta t}\right|\right|\leq C,
    \end{equation}
    (where $U^{n}$ are the numerical density of each species $n_p^n$, the total velocity $u^n$ and the total energy $E^n$ of the mixture computed at time $t^n$) the scheme asymptotically becomes a first order in time approximation of the compressible Navier-Stokes for gas mixtures given by
    \begin{equation}\label{eqNS}
        \begin{cases}
            \displaystyle \frac{n_p^{n+1}-n_p^n}{\Delta t} + \nabla_x\cdot (n_p^n u^n) = -\varepsilon \nabla_x\cdot\left(J_p^{n-1}\right), \quad \quad \forall p,\\
            \displaystyle \frac{\rho^{n+1}u^{n+1}-\rho^n u^n}{\Delta t} + \nabla_x\cdot (\rho^n u^n \otimes u^n) + \nabla_x(n^n T^n) = - \varepsilon \nabla_x\cdot\left(\rho^{n-1} \Theta^{n-1}\right), \\
            \displaystyle \frac{E^{n+1}-E^n}{\Delta t} + \nabla_x\cdot\left(\left[E^n + n^n T^n\right]u^n\right) = - \varepsilon \nabla_x\cdot\left(\mathbb{Q}^{n-1} + \rho^{n-1} \Theta^{n-1} u^{n-1}\right),
        \end{cases}
    \end{equation}
\end{enumerate}
where $\Theta^n$ is given by equation \eqref{eq6000}, while $J_p^n$ and $\mathbb{Q}^n$ at time $t^n$ are
\begin{equation}
    \begin{cases}
        \displaystyle \mathbb{Q}^n = \sum_{p=1}^P \int m_p f_p^n \frac{|v-u^n|^2}{2}(v-u^n) \  dv,\\
        \displaystyle J_p = \int vf_p^n \  dv,
    \end{cases}
\end{equation}
and they will be computed in this section providing their analytical expression.\\

\noindent\emph{Proof} To prove the first property, let's note that $f_p^{n+1}$ is a linear combination of $f^n_p$ and $G_{p}^{n+1}$ thus we get the first assertion.\\
To prove the second property let us consider a generic initial data and $f^n_p$ for $n\geq 1$ and let us compute the asymptotic limit of $\Sigma^n$ when $\varepsilon$ goes to zero in equation \eqref{eq1}. It yields
\begin{equation}
    \Sigma^n = T^{n}n^{n}\text{I} + \rho^{n} u^n \otimes u^n,
\end{equation}
and using \eqref{eq2}
\begin{equation}
    \rho^n \Theta^n = T^n  n^n \text{I},
\end{equation}
and from \eqref{eq3}
\begin{equation}
    \tilde{\tau}_p^n =  T^n \frac{1}{m_p} \text{I}.
\end{equation}
So in the asymptotic limit $\varepsilon\to0$ the quantities $G_p^n$ become isotropic, which means that they converge to the corresponding equilibrium Maxwellians $\mathcal{M}^{n_p^n, u^n, T^n}_p$. So the solution at zeroth order is obtained by taking $f^n_p=\mathcal{M}^{n_p^n, u^n, T^n}_p$ in the conservation laws \eqref{eq4}
\begin{equation}
    \begin{cases}
        \displaystyle n_p^{n+1} = \int \left(\mathcal{M}_p^{n_p^n, u^n, T^n} - \Delta t \  v\cdot \nabla_x \mathcal{M}_p^{n_p^n, u^n, T^n}\right)\  dv,\\
        \displaystyle \rho^{n+1}=\sum_{p=1}^P \int m_p \left(\mathcal{M}_p^{n_p^n, u^n, T^n} - \Delta t \  v\cdot \nabla_x \mathcal{M}_p^{n_p^n, u^n, T^n}\right)\  dv=\sum_p m_p n_p^{n+1},\\
        \displaystyle \rho^{n+1} u^{n+1} = \sum_{p=1}^P \int m_p v \left(\mathcal{M}_p^{n_p^n, u^n, T^n} - \Delta t \  v\cdot \nabla_x \mathcal{M}_p^{n_p^n, u^n, T^n}\right)\  dv,\\

        \displaystyle E^{n+1} = \sum_{p=1}^P \int \frac{m_p}{2}|v|^2 \left(\mathcal{M}_p^{n_p^n, u^n, T^n} - \Delta t \  v\cdot \nabla_x \mathcal{M}_p^{n_p^n, u^n, T^n}\right) dv,
    \end{cases}
\end{equation}

from which one gets the following approximation of the multi-species Euler equation
\begin{equation}\label{eqEuler}
    \begin{cases}
        \displaystyle  \frac{n_p^{n+1}-n_p^n}{\Delta t} + \nabla_x\cdot (n_p^{n} u^{n})= 0, \quad \quad \forall p\\
        
        \displaystyle \frac{\rho^{n+1}u^{n+1} - \rho^{n}u^n}{\Delta t} + \nabla_x \cdot (\rho^{n} u^{n}\otimes u^{n})+\nabla_x(n^{n} T^{n})=0,\\
        \displaystyle \frac{E^{n+1}-E^n}{\Delta t} + \nabla_x \cdot \left[\left(\frac{1}{2}\rho^n |u^n|^2 + \frac{5}{2}n^n T^n\right)u^n\right]=0.
    \end{cases}
\end{equation}
This concludes the proof of the second assertion.\\
Now let's prove the third property regarding the asymptotically convergence towards the Navier-Stokes equation, using the Chapman-Enskog expansion, which consists of expanding each distribution function $f_p^n$ as
\begin{equation}
    f_p^n = \mathcal{M}^n_p + \varepsilon h_p^n,
\end{equation}
frow which one gets
\begin{equation}\label{eq5000}
    \begin{cases}
        \displaystyle\int h_p^n\  dv = 0 \quad \quad \forall p,\\
        \displaystyle \sum_{p=1}^P \int m_p v h_p^n\  dv = \sum_{p=1}^P \int m_p |v|^2 h_p\  dv =0.
    \end{cases}
\end{equation}
Similarly we also expand
\begin{equation}\label{eq3000}
    \Theta^n = \frac{1}{\rho^n}\sum_{p=1}^P \int m_p (v-u^n)\otimes(v-u^n) \  f_p^n \  dv = \frac{T^n n^n}{\rho^n}\text{I} + \varepsilon \Theta_1^n,
\end{equation}
and
\begin{equation}
    \displaystyle \mathbb{Q}^n = \sum_{p=1}^P\int  m_p\frac{|v-u^n|^2}{2}(v-u^n)f_p^n \  dv = 0 + \varepsilon \mathbb{Q}_1^n,
\end{equation}
where
\begin{equation}\label{eqBoth}
    \begin{cases}
        \displaystyle \rho^n \Theta^n_1 = \sum_{p=1}^P \int m_p (v-u^n)\otimes (v-u^n) h_p^n \  dv,\\
        \displaystyle \mathbb{Q}_1^n = \sum_{p=1}^P \int m_p \frac{|v-u^n|^2}{2}(v-u^n) h_p^n \  dv,\\
    \end{cases}
\end{equation}
where $\Theta^n_1$ is traceless, because of equation \eqref{eq5000}.\\
Let us note that
\begin{equation}
    J_{p,1} = \int h_pv\  dv=\int h_p(v-u)\  dv,
\end{equation}
where in the last equality we use equation \eqref{eq5000}.
We can use these quantities inside the conservations laws \eqref{eq4}, obtaining 
\begin{equation}\label{eq3010}
    \begin{cases}
        \displaystyle \frac{n_p^{n+1}-n_p^n}{\Delta t} + \nabla_x\cdot (n_p^n u^n) = -\varepsilon \nabla_x\cdot\left(J_p^{n}\right), \quad \quad \forall p,\\
            \displaystyle \frac{\rho^{n+1}u^{n+1}-\rho^n u^n}{\Delta t} + \nabla_x\cdot (\rho^n u^n \otimes u^n) + \nabla_x(n^n T^n) = - \varepsilon \nabla_x\cdot\left(\rho^{n} \Theta^{n}_1\right), \\
            \displaystyle \frac{E^{n+1}-E^n}{\Delta t} + \nabla_x \cdot \left[\left(\frac{1}{2}\rho^n |u^n|^2 + \frac{5}{2}n^n T^n\right)u^n\right] = - \varepsilon \nabla_x\cdot\left(\mathbb{Q}^{n}_1 + \rho^{n} \Theta^{n}_1 u^{n}\right),
    \end{cases}
\end{equation}
In the same way, we also expand each anisotropic Gaussian
\begin{equation}
    G_p^n = \mathcal{M}_p^n + \varepsilon g_p^n.
\end{equation}
From the definition of $\tau$ \eqref{eq3} and \eqref{eq3000} one gets that
\begin{equation}
    \tau^n = T^n \frac{n^n}{\rho^n} \text{I}+\nu \varepsilon \Theta_1^n.
\end{equation}
Using equation \eqref{eq3}
\begin{equation}
    \tilde{\tau}_p^n = \frac{T^n}{m_p}\text{I} + \nu\varepsilon \Theta_1^n \frac{\rho^n}{m_p n^n}, 
\end{equation}
whose determinant is equal to $(\frac{T^n}{m_p})^{D_V} + \mathcal{O}(\varepsilon^2)$ \cite{fil-jin-2011}, and the inverse is \cite{fil-jin-2011}
\begin{equation}
    [\tilde{\tau}_p^n]^{-1} = \frac{m_p}{T^n}
    \left[\text{I} - \frac{\nu \varepsilon \rho^n}{n^n T^n }\Theta_1^n\right]  + \mathcal{O}(\varepsilon^2).
\end{equation}
The quantity $g_p^n$ can be computed as
\begin{equation}\label{eq3001}
    g_p^n = \frac{G_p^n - \mathcal{M}_p^n}{\varepsilon} = \frac{1}{2}\frac{m_p}{n^n (T^n)^2}\nu \rho^n (v-u^n)^\intercal \Theta_1^n (v-u^n) \mathcal{M}_p^n.
\end{equation}
Considering the terms of order zero into the scheme \eqref{eq5}, one get an expression for $h_p^n$
\begin{equation} \label{eq3003}
    h_p^n = g_p^n - \frac{\mathcal{M}_p^n - \mathcal{M}_p^{n-1}}{\lambda^n \Delta t} - \frac{v}{\lambda^n}\cdot \nabla_x \mathcal{M}_p^{n-1} + \mathcal{O}(\varepsilon).
\end{equation}
From \cite{bis-gro-mar-2021} one gets
\begin{equation}\label{eq3002}
\begin{aligned}
    &\left(\partial_t + v \nabla_x\right)\mathcal{M}_p=\\
    &= \mathcal{M}_p
    \left[ m_p \  c\cdot \left(\frac{\nabla_x n_p}{\rho_p}-\frac{\nabla_x n}{\rho}\right) + \frac{m_p}{T}\left(c\otimes c - \frac{1}{3}|c|^2 \text{I}\right):\nabla_x u + \left(\frac{m_p}{2T}\left|c\right|^2 - \frac{3}{2} - \frac{m_p n}{\rho}\right)c\cdot \frac{\nabla_x T}{T}\right],
\end{aligned}
\end{equation}
where $c = v - u$.\\

Combining equations \eqref{eq3001} and \eqref{eq3002} into equation \eqref{eq3003}, and using assumption \eqref{eqAssumption}, it is possible to get an expression for $h_p^n$
\begin{equation}
\begin{aligned}
    h_p^n =&\\
    =& - \frac{1}{\lambda^n} \mathcal{M}_p^{n-1}
    \Biggl[ m_p \  c^{n-1}\cdot \left(\frac{\nabla_x n_p^{n-1}}{\rho_p^{n-1}}-\frac{\nabla_x n^{n-1}}{\rho^{n-1}}\right) + \frac{m_p}{T^{n-1}}\left(c^{n-1}\otimes c^{n-1} - \frac{1}{3}|c^{n-1}|^2 \text{I}\right):\nabla_x u^{n-1} \\
    &+\left(\frac{m_p}{2T^{n-1}}|c^{n-1}|^2 - \frac{3}{2} - \frac{m_p n^{n-1}}{\rho^{n-1}}\right)c^{n-1}\cdot \frac{\nabla_x T^{n-1}}{T^{n-1}}\Biggr]  \\
    &+ \frac{1}{2}\frac{m_p}{n^n (T^n)^2}\nu \rho^n (c^n)^\intercal \Theta_1^n (c^n) \mathcal{M}_p^n + \mathcal{O}(\Delta t) + \mathcal{O}(\varepsilon).
\end{aligned}
\end{equation}
It's now possible to compute $\Theta_1^n$ and $\mathbb{Q}_1^n$ at the first order of the Chapman-Enskog expansion, by substituting this last expression for $h_p^n$ into equation \eqref{eqBoth}.
\begin{equation}
    \rho^n \Theta_1^n = \sum_p \int m_p (v-u^n)\otimes(v-u^n)h_p^n \  dv = -\frac{n^{n-1}T^{n-1}}{\lambda^n}\sigma(u^{n-1}) + \nu \Theta_1^n \rho^n  + \mathcal{O}(\Delta t) + \mathcal{O}(\varepsilon),
\end{equation}
where $\sigma(u) = \nabla_x u + \left(\nabla_x u\right)^\intercal - \frac{2}{3}\left(\operatorname{div}_x u\right) \text{I}$.\\
Rearranging the terms one gets
\begin{equation}
    \rho^n\Theta_1^n = -\frac{n^{n-1}T^{n-1}}{(1-\nu)\lambda^n}\sigma(u^{n-1}) = -\mu \sigma(u^{n-1})  + \mathcal{O}(\Delta t) + \mathcal{O}(\varepsilon),
\end{equation}
where we used the fact that the viscosity \cite{bru-2015} is given by
\begin{equation}\label{eqViscosity}
    \mu = \frac{nT}{(1-\nu)\lambda}.
\end{equation}
Similarly
\begin{equation}
    \mathbb{Q}_1^n = \sum_{p=1}^P \int m_p \frac{|v-u^n|^2}{2}(v-u^n)h_p^n \  dv = \frac{5}{2}T^{n-1}\left(\sum_p J_{p,1}^n\right) - k \nabla_x T^{n-1}  + \mathcal{O}(\Delta t) + \mathcal{O}(\varepsilon),
\end{equation}
where the thermal conductivity $k$, consistently with \cite{bru-2015}, is
\begin{equation}
    k = \frac{5}{2}\frac{1}{\lambda^n}\left(\sum_p\frac{n_p^{n-1}}{m_p}\right)T^{n-1},
\end{equation}
and 
\begin{equation}
    J_{p,1}^n = - \frac{n_p^{n-1}}{\lambda^{n}} \left(\frac{\nabla_x n_p^{n-1}}{\rho_p^{n-1}} - \frac{\nabla_x n^{n-1}}{\rho^{n-1}}\right) T^{n-1} - \frac{n_p^{n-1}}{\lambda^{n}} \left(\frac{1}{m_p} - \frac{n^{n-1}}{\rho^{n-1}}\right)\nabla_x T^{n-1}  + \mathcal{O}(\Delta t) + \mathcal{O}(\varepsilon).
\end{equation}
The Prandtl Number is related to the coefficient $\nu$ by \cite{fil-jin-2011, bru-2015}
\begin{equation}
    \text{Pr} = \frac{5}{2}\frac{\mu}{k} = \frac{n^{n-1}}{1-\nu}\frac{1}{\sum_p\frac{n_p^{n-1}}{m_p}}.
\end{equation}

\subsection{High-order IMEX shemes}
The IMEX scheme showed in the previous section is first order in time. However it can be extended to an arbitrary higher temporal order using, for example, the IMEX DIRK (Diagonally Implicit Runge Kutta) schemes proposed by \cite{par-rus-2001}, that we briefly review in the following (much more details on the topic can be found in \cite{bos-par-rus-2024}).\\
The approach presented in \cite{par-rus-2001} is adoptable for any stiff systems of differential equation in the form
\begin{equation}\label{eq1001}
    y'=q(y) + \frac{1}{\varepsilon} g(y),
\end{equation}
where $y=y(t)\in \RR^N$, $q, g: \RR^N\to \RR^N$. In our case $y$ is one of the distribution functions of the mixture, when considering only the temporal dependence, $q(y)$ is the transport part of the Boltzmann equation, and $g(y)$ is the ES-BGK operator. In order to extend this procedure to the gas mixture, then we have to apply this scheme to each distribution function.\\

An Implicit-Explicit IMEX Runge Kutta scheme for system \eqref{eq1001} is of the form
\begin{equation}
    \begin{cases}
        \displaystyle Y_i = y_n + \Delta t \sum_{j=1}^{i-1} \tilde{a}_{ij}\  q(t_n + \tilde{c}_j \Delta t, Y_j) + \Delta t \sum_{j=1}^\Lambda a_{ij}\frac{1}{\varepsilon} g(t_n + c_j \Delta t, Y_j),\\
        \displaystyle y_{n+1}=y_n + \Delta t \sum_{i=1}^\Lambda \tilde{\omega}_{i} \  q(t_n+\tilde{c}_i \Delta t, Y_i)+\Delta t \sum_{i=1}^\Lambda \omega_i \frac{1}{\varepsilon} g(t_n + c_i \Delta t, Y_i).
    \end{cases}                    
\end{equation}
The matrices $\tilde{A}=(\tilde{a}_{ij})$, $\tilde{a}_{ij}=0$ for $j\geq i$ and $A=(a_{ij})$ are $\Lambda \times \Lambda$ matrices such that the scheme is explicit in $q$, and implicit in $g$. An IMEX Runge-Kutta scheme is characterized by these two matrices and the coefficient vectors $\tilde{c}=(\tilde{c}_1, ..., \tilde{c}_\Lambda)^\intercal$, $\tilde{\omega}=(\tilde{\omega}_1, ..., \tilde{\omega}_\Lambda)^\intercal$, $c=(c_1, ..., c_\Lambda)^\intercal$, $\omega = (\omega_1, ..., \omega_\Lambda)^\intercal$, and it can be represented by a double \emph{tableau} in the Butcher notation,
\begin{center}
    \begin{tabular}{c|c}
        $\tilde{c}$ & $\tilde{A}$\\                
        \hline
        \\
                    & $\tilde{\omega}^\intercal$\\
    \end{tabular}, \  \   
    \begin{tabular}{c|c}
        $c$ & $A$\\                
        \hline
        \\
            & $\omega^\intercal$\\
    \end{tabular}.
\end{center}
In our code we implement the third order time IMEX DIRK scheme ARS(2,3,3) \cite{par-rus-2001}, which is characterized by the following Butcher tableaus
\begin{center}
    \begin{tabular}{c|c c c}
        $0$        & $0$        & $0$         & $0$\\   
        $\gamma$   & $\gamma$   & $0$         & $0$\\ 
        $1-\gamma$ & $\gamma-1$ & $2-2\gamma$ & $0$\\
        \hline
                   & $0$          & $1/2$       & $1/2$\\
    \end{tabular}, \  \  
    \begin{tabular}{c|c c c}
        $0$        & $0$ & $0$         & $0$\\   
        $\gamma$   & $0$ & $\gamma$    & $0$\\ 
        $1-\gamma$ & $0$ & $1-2\gamma$ & $\gamma$\\
        \hline
                   & $0$          & $1/2$       & $1/2$\\
    \end{tabular},
\end{center}
where $\gamma = (3+\sqrt3)/6$.

\subsection{Conservation of macroscopic quantities}
It is well known that the discretization of the velocity space and the introduction of velocity grid introduce a loss of conservation \cite{dim-par-2014}. A natural solution to the problem is to consider a velocity support big enough such that the distribution function is smooth and that the energy outside the domain in negligible. However, it may be considerable to construct a solver where conservation properties are exactly maintained during the simulation.\\
Several approaches have been proposed in literature, exploiting an Entropy minimization (like in \cite{mie-2000}), or a $L^2$-minimization (like in \cite{gam-tha-2009, bos-cho-rus-2022}). The \emph{Entropy approach} allows the construction of a conservative discrete Maxwellian, while $L^2-$minimization can be applied to more general distribution functions. In this paper we implement, describe and extend to a multi-species context, the weighted $L^2-$minimization approach proposed in \cite{bos-cho-rus-2022}. In particular, we generalize and extended this method to a multi-species context applying it to each anisotropic Maxwellian $G_p$ computed with our scheme. We will present the method considering an arbitrary quantity of species, and in a spatially homogenous settings. Generalizing to a spatially non homogenous settings is straightforward: it is sufficient to apply the following strategy cell-by-cell in the space domain. In the following we will consider a space of dimension 3 in velocity, but of course it can be generalized to an arbitrary amounts of dimensions in the velocity space.\\
Let $\tilde{G}_p$ be the initial guess for the anisotropic Maxwellian of the $p-$th species obtained with the previous AP-IMEX scheme. Let us consider a different generic weight function $1/h_p$ for each species (we will specify in the following the choice of the weight functions), then we look for the target anisotropic Maxwellians $G_p$ exactly preserving the moments.
Let $M$ be the total number of discretization points of the velocity space such that
\begin{equation}
    \tilde{G}_p = \begin{pmatrix}
           \tilde{G}_p^1\\
           \tilde{G}_p^2\\
           ...\\
           \tilde{G}_p^M
          \end{pmatrix}. 
\end{equation}
Let $Z_p$ be the $(D_V + 2)\times M$ (in this case $D_V=3$, number of dimensions of the velocity space) integration matrix
\begin{equation}
\begin{aligned}
    Z_p &= \begin{pmatrix}
           w_1 (h_p)_1                & w_2 (h_p)_2                 & ...& w_M (h_p)_M\\
           m_p v^x_1 w_1 (h_p)_1      & m_p v^x_2 w_2 (h_p)_2       & ...& m_p v^x_M w_M (h_p)_M\\
           m_p v^y_1 w_1 (h_p)_1      & m_p v^y_2 w_2 (h_p)_2       & ...& m_p v^y_M w_M (h_p)_M\\
           m_p v^z_1 w_1 (h_p)_1      & m_p v^z_2 w_2 (h_p)_2       & ...& m_p v^z_M w_M (h_p)_M\\
           0.5 m_p |v_1|^2w_1 (h_p)_1 & 0.5 m_p |v_2|^2 w_2 (h_p)_2 & ...& 0.5 m_p |v_M|^2 w_M (h_p)_M\\
          \end{pmatrix},
\end{aligned}    
\end{equation}
where $w_i$ are the integration weights.\\
The moments of the $p-$th species are the following
\begin{equation}
    a_p = \begin{pmatrix}
         n_p\\
         \rho_p u_p^x\\
         \rho_p u_p^y\\
         \rho_p u_p^z\\
         E_p
        \end{pmatrix},
\end{equation}
where $n_p$ is the numerical density of the $p$-th species, $\rho_p$ is its mass density, $u_p^x, u_p^y, u_p^z$ are respectively, the total velocities of the $p$-th species in the three directions and $E_p$ is its energy.\\
The goal is to find $G_p$ that minimizes the \emph{weighted} $L^2$ distance from $\tilde{G}_p \circ \frac{1}{h_p}$ under the constraint that $Z_pG_p - a_p=0$, where the symbol $\circ$ is the Hadamard product.\\
Similarly to \cite{bos-cho-rus-2022}, we introduce the following Lagrangian for each species
\begin{equation}
        \displaystyle L_p = \left|\left|\tilde{G}_p\circ\frac{1}{h_p} - G_p\right|\right|^2_2 + \lambda_p^\intercal(Z_pG_p  - a_p).
\end{equation}
From the stationary points of $L_p$ one gets an expression for $\lambda_p$ from which \cite{bos-cho-rus-2022}
\begin{equation} 
        \displaystyle G_p \circ h_p = \tilde{G}_p + Z_p^\intercal (Z_p Z_p^\intercal)^{-1} \left(a_p - Z_p \left( \tilde{G}_p \circ \frac{1}{h_p}\right)\right) \circ h_p.
\end{equation}
Then in order to impose the conservation of moments we applied this approach to each anisotropic Maxwellian $\tilde{G}_p$ computed using our scheme. In particular, at each step of the Runge-Kutta scheme, we compute the anisotropic Maxwellians $\tilde{G}_p$. Then we use this $L^2-$minimization strategy, to estimate the anisotropic Maxwellians (the quantity $G_p$ in this last approach) having exactly the desired moments used to construct $\tilde{G}_p$ using our scheme. Similarly to \cite{bos-cho-rus-2022}, we impose $h_p$ to be equal to the Maxwellian $\mathcal{M}_p$ having the moments used to construct $\tilde{G}_p$.

\section{Numerical tests}\label{secNumTest}
In the following, we present several simulations conducted to numerically test and validate our numerical scheme. All simulations have been conducted on a non-dedicated computer equipped with a \emph{12th Gen Intel(R) Core(TM) i7-12700H} CPU and 16 GB of RAM, except the \hyperref[Test2]{\textbf{Test 2}} which has been executed on a dedicated computer equipped with two \emph{Intel(R) Xeon(R) Gold 5418Y} CPUs and 128 GB of RAM, kindly provided by the Mathematics Department of the University of Ferrara. In order to speed-up the code execution a multithread approach has been implemented for the non homogeneous simulations 2D in space: each thread operates on a subdomain of the entire spatial domain.
For the spatially non-homogeneous simulations the variables in the transport term have been reconstructed using a CWENO 3 scheme \cite{Lev-Pup-Rus-1999, pup-sem-vis-2023}, while in time we use the IMEX scheme presented before extended to the third order using the ARS(2,3,3) scheme \cite{par-rus-2001}. For the sake of simplicity, in all simulations we considered $\nu=-0.5$ and $\lambda=1$, except where clearly indicated.\\

In Section \ref{validation} we perform some numerical experiments to validate the scheme and our numerical implementation. In Section \ref{applications} we show some interesting numerical applications for fluid simulations 1D and 2D in space, like the classical SOD shock tube benchmark test, or some intriguing phenomena, like the Kelvin-Helmholtz instability and the flow of a fluid around a cylinder. 

\subsection{Validation}\label{validation}
In this section we perform some numerical experiments to validate the scheme and our numerical implementation. We perform a first \hyperref[Test00]{\textbf{test}} in order to verify the moments relaxation of two species toward those of equilibrium. In a second numerical \hyperref[Test01]{\textbf{experiment}} we study the convergence rate in time and both in space and time of our code.

\subsubsection{Homogeneous relaxation test}\label{Test00}
\begin{figure}
    \centering
    \subfigure[Temperature Evolution; Mass Ratio = 1]{
        \includegraphics[scale=0.54]{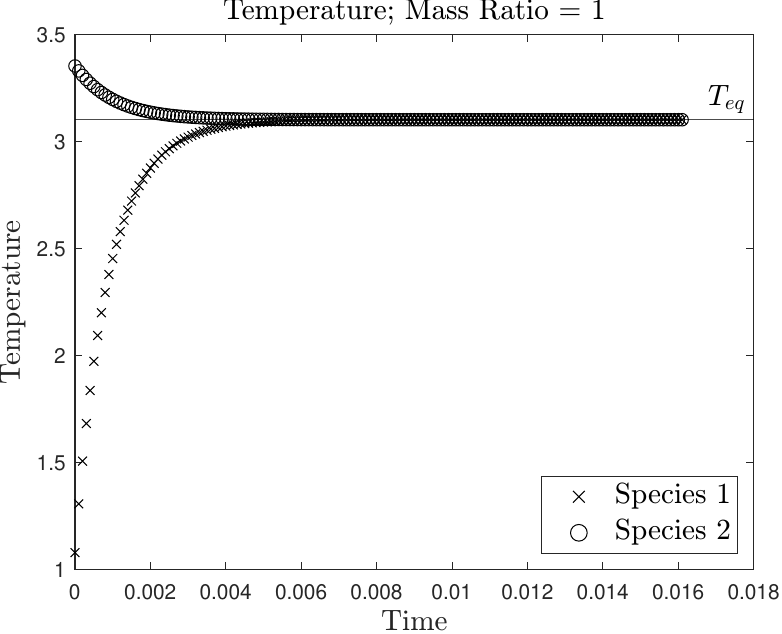}
        \label{fig:enter-label1} 
    }
    \subfigure[Velocity Evolution; Mass Ratio = 1]{
        \includegraphics[scale=0.54]{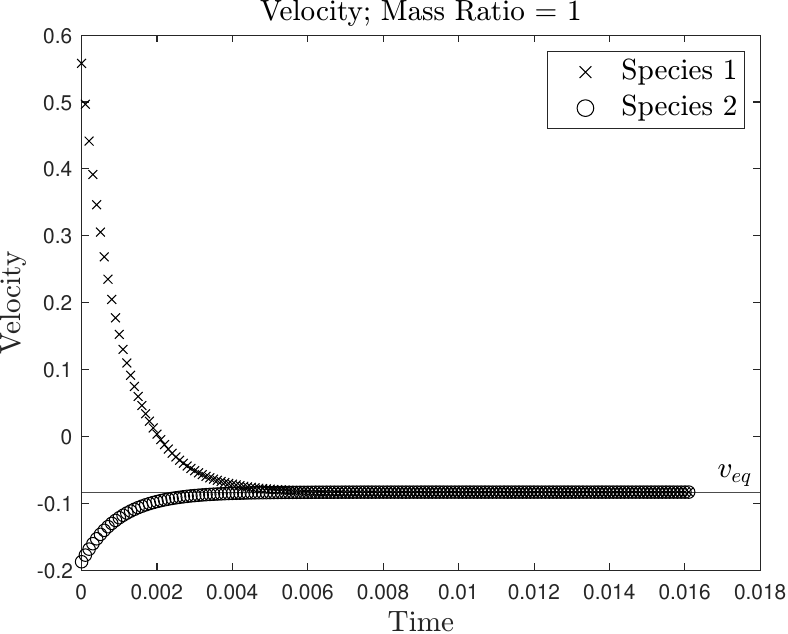}
        \label{fig:enter-label2} 
    }    
    \vspace{0.5cm}
    \subfigure[Temperature Evolution; Mass Ratio = 100]{
        \includegraphics[scale=0.54]{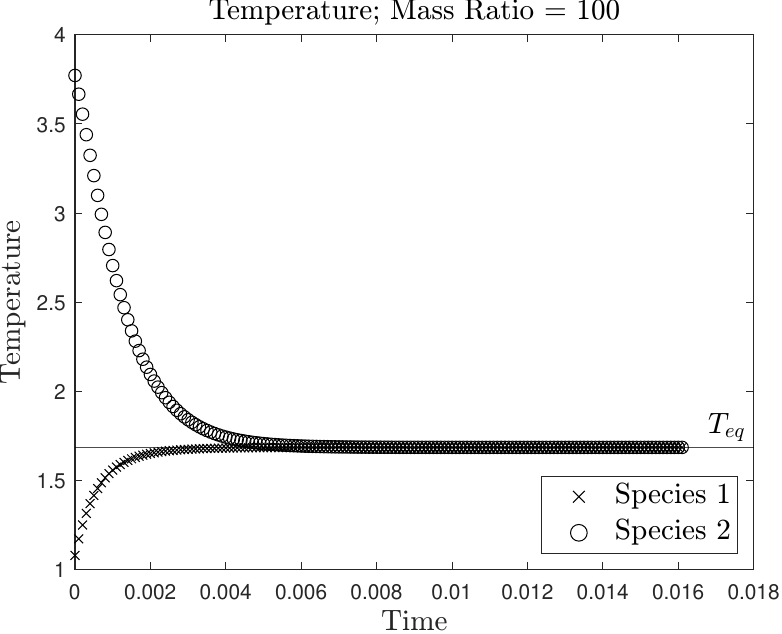}
        \label{fig:enter-label3} 
    }
    \subfigure[Velocity Evolution; Mass Ratio = 100]{
        \includegraphics[scale=0.54]{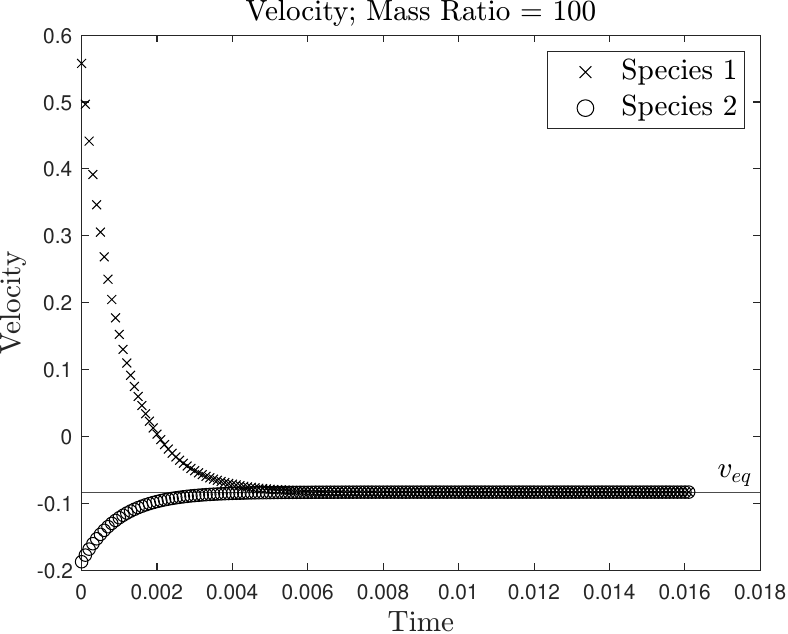}
        \label{fig:enter-label4} 
    }
    \caption{Moments evolution for the homogeneous relaxation \hyperref[Test00]{\textbf{test}}. In the first row it is reported the evolution of the temperature and the velocity of each species for a mass ratio equal to 1. The same quantities, but for a mass ratio equal to 100, are reported in the second row. The solid line in each figure represents the asymptotic value associated to the corresponding equilibrium Maxwellian given by equation \eqref{eq7}, while the crosses and the open circles represents, respectively, the first and the second species.}    
    \label{fig-Test0}
\end{figure}

The aim of this spatially homogeneous test 1D in velocity is to check the correct convergence of the temperature of each species toward the equilibrium one. We have considered two species of masses $m_1$ and $m_2$ described by the following distribution functions $f_1$ and $f_2$ both initialized as sum of two Maxwellians
\begin{equation}
    \begin{cases}
        f_1(0, v) = \mathcal{M}_1^{n_a, v_a, T_a} + \mathcal{M}_1^{0.3 n_a, 1.5 v_a, 1.3 T_a},\\
        f_2(0, v) = \mathcal{M}_2^{n_b, v_b, 0.7T_b} + \mathcal{M}_2^{3n_b, 0.5 v_b, 2 T_b},\\
        n_a = 1/m_1, \quad T_a=1, \quad v_a=0.5,\\
        n_b = 2/m_2, \quad T_b=2, \quad v_b=-0.3.\\
    \end{cases}    
\end{equation}
It is important to remark that the macroscopic quantities $n_a$, $n_b$, $v_a$, $v_b$, $T_a$ and $T_b$ do not correspond to the macroscopic moments of the two species $n_1$, $n_2$, $v_1$, $v_2$, $T_1$ and $T_2$, since both species are sum of two different Maxwellians. Also, the Maxwellians $\mathcal{M}_p$ are computed considering the mass of the $p$-th species. Thus in the expression of the Maxwellians $\mathcal{M}_1$, the mass $m_1$ is employed; similarly, in the expression of the Maxwellians $\mathcal{M}_2$, the mass $m_2$ is employed.\\ 
The simulations have been conducted in two different configurations of mass ratio ($MR$): configuration 1 with $m_2/m_1=1$ and configuration 2 with $m_2/m_1=100$. In all simulations and both configurations we considered $m_1=1$, $\varepsilon=10^{-3}$ (Knudsen Number), $\Delta t=10^{-4}$ (time step). The small value of $\Delta t$ is not due to any stability condition, but only in order to better graphically visualize the convergence of the parameters of each species towards the one of the equilibrium Maxwellian. The final time of the integration is equal to $16\varepsilon$. We considered a velocity box $[-20, 20]$ for both configurations, discretized with a total of 256 points when $MR=100$ and 32 points when $MR=1$. In Figure \ref{fig-Test0}, the temperature and velocity evolution in time of both species are plotted, and it is evident that they are converging toward the corresponding asymptotic value.\\
The time integration in this simulation has been performed using the third-order IMEX strategy presented above, without considering the transport term in the Boltzmann equation.

\subsubsection{Convergence test}\label{Test01}
\begin{figure}
    \centering
    \includegraphics[scale=0.8]{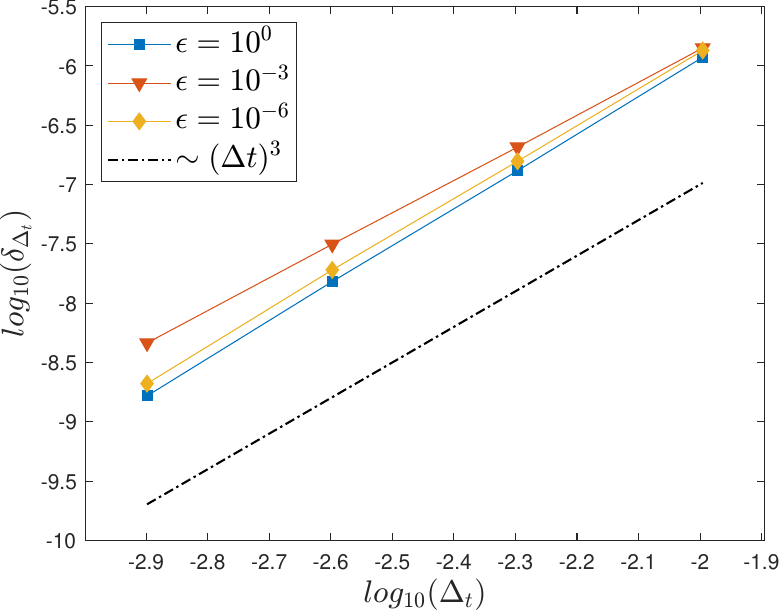}
    \caption{Convergence rate in time for different values of $\varepsilon$ for the time convergence \hyperref[Test01]{\textbf{test}}.} 
    \label{fig_convTest}
\end{figure}

\paragraph{Time convergence test}
The aim of this test is to check the third-order temporal convergence of the implemented ARS(2,3,3) IMEX scheme. We considered two gases with a given initial condition 1D in space and 1D in velocity, and we evolved the solution using the IMEX scheme presented above, using different timesteps $(\Delta t)_i$ up to a time $T_f=1$ (where $(\Delta t)_{i+1}=0.5 (\Delta t)_i$). Then we evaluated the error $\delta_{(\Delta t)_i}$ associated to each $(\Delta t)_i$ at the final integration time $T_f$ as following
\begin{equation}
    \delta_{(\Delta t)_i} = \left|\left| \rho_{(\Delta t)_i}(T_f) - \rho_{(\Delta t)_{rif}}(T_f) \right|\right|_2,
\end{equation}
where $\rho_{(\Delta t)_{rif}}$ is the solution computed using the timestep $(\Delta t)_{rif}=0.5 (\Delta t)_{min}$, and $(\Delta t)_{min}$ is the smallest timestep for which we compute the error convergence.\\
We considered two gases of masses $m_2=2 m_1=2$, initialised as
\begin{equation}
    \begin{cases}
        f_1 = \mathcal{M}_1^{n_1(x), v(x), T(x)}\\
        f_2 = \mathcal{M}_2^{n_2(x), v(x), T(x)},
    \end{cases}
\end{equation}
where $v(x)=0$, $T(x)=0.3 \left(1+\exp{(-3|x-L_x/2|^2)}\right)$, $n_i(x)=\rho(x)/m_i$ and\\ $\rho(x)=1+0.1\exp{(-3|x-L_x/2|^2)}$, where $L_x=10$, and $x\in [0, L_x]$ discretized using 100 cells. The velocity space $[-10, 10]$ is discretized using 40 points, and we used periodic boundary conditions in the space domain. The convergence curves for different values of $\varepsilon$ are reported in Figure \ref{fig_convTest}. They clearly show the third order of convergence of the scheme.

\begin{figure}
    \centering
    \includegraphics[scale=0.8]{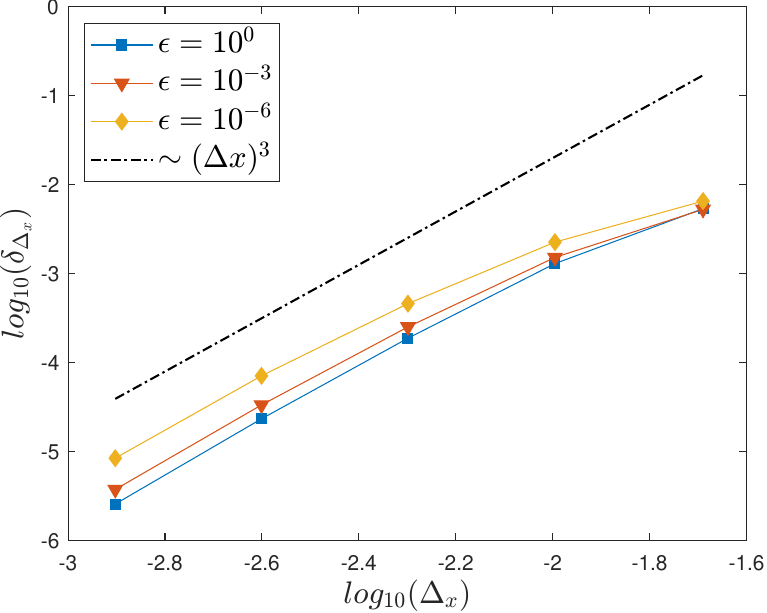}
    \caption{Convergence rate in space and time for different values of $\varepsilon$ for the space and time convergence \hyperref[Test01]{\textbf{test}}.} 
    \label{fig_convTest_space_time}
\end{figure}

\paragraph{Space and time convergence test}
Lastly we have also checked the global third-order accuracy of the scheme in space and in time.
We considered two gases with a given initial condition 1D in space and 1D in velocity, and we evolved the solution using the IMEX scheme presented above, using different spatial discretization steps $(\Delta x)_i$ (obtained by doubling the number of discretization cells $N_{i+1}=2 N_i$). The timestep associated to each run is $(\Delta t)_i = (\Delta x)_i / L_v$ (where $L_v=10$ is the maximum value of the velocity box $[-10, 10]$ discretized using 40 points). Then we evaluated the error $\delta_{(\Delta x)_i}$ associated to each $(\Delta x)_i$ at the final integration time $T_f=0.2$ as following
\begin{equation}
    \delta_{(\Delta x)_i} = \left|\left| \rho_{(\Delta x)_i}(T_f) - \rho_{(\Delta x)_{rif}}(T_f) \right|\right|_2,
\end{equation}
where $\rho_{(\Delta x)_{rif}}$ is the solution computed using the timestep $(\Delta x)_{rif}= 0.5(\Delta x)_{min}$, and $(\Delta x)_{min}$ is the smallest step for which we compute the error convergence.\\
We considered two gases of masses $m_2=2 m_1=2$, initialised as
\begin{equation}
    \begin{cases}
        f_1 = \mathcal{M}_1^{n_1(x), v(x), T(x)}\\
        f_2 = \mathcal{M}_2^{n_2(x), v(x), T(x)},
    \end{cases}
\end{equation}
where $v(x)=0$, $T(x)=0.3 \left(1+0.1\exp{(-\Theta|x-L_x/2|^2)}\right)$, $n_i(x)=\rho(x)/m_i$ and\\ $\rho(x)=1+0.1\exp{(-\Theta|x-L_x/2|^2)}$, where $\Theta=1000$, $L_x=1$, and $x\in [0, L_x]$. The first run $i=1$ is the one obtained using 50 cells to discretize the physical space. We used periodic boundary conditions in the space domain. In Figure \ref{fig_convTest_space_time} the convergence curves for different values of $\varepsilon$ are reported. They clearly show the global third order of convergence of the scheme.

\subsection{Numerical applications}\label{applications}
We show some interesting numerical applications for fluid simulations 1D and 2D in space, like the classical SOD shock tube benchmark test, or some intriguing phenomena, like the Kelvin-Helmholtz instability and the flow over a cylinder. Lastly, we perform a comparison between BGK, ES-BGK and full Boltzmann operator. This last test shows that the ES-BGK performs better than the BGK operator in the rarefied regime. At the end of this section we also discuss the big speed-up offered by such kind of IMEX scheme compared to fully explicit integrators.
\subsubsection{Test 1: A Sod shock tube problem}\label{Test1}
\begin{figure}
    \centering
    \subfigure[Density]{
        \includegraphics[scale=0.56]{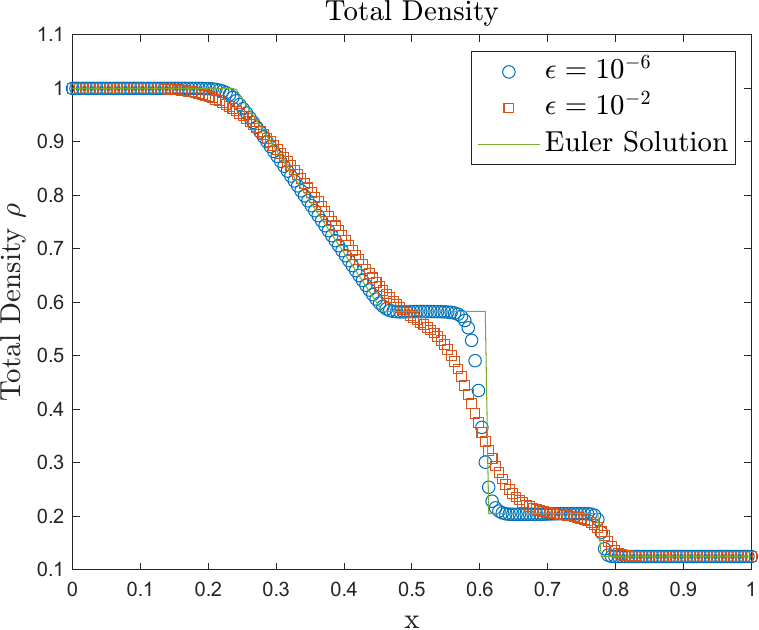}
    }
    \subfigure[Temperature]{
        \includegraphics[scale=0.56]{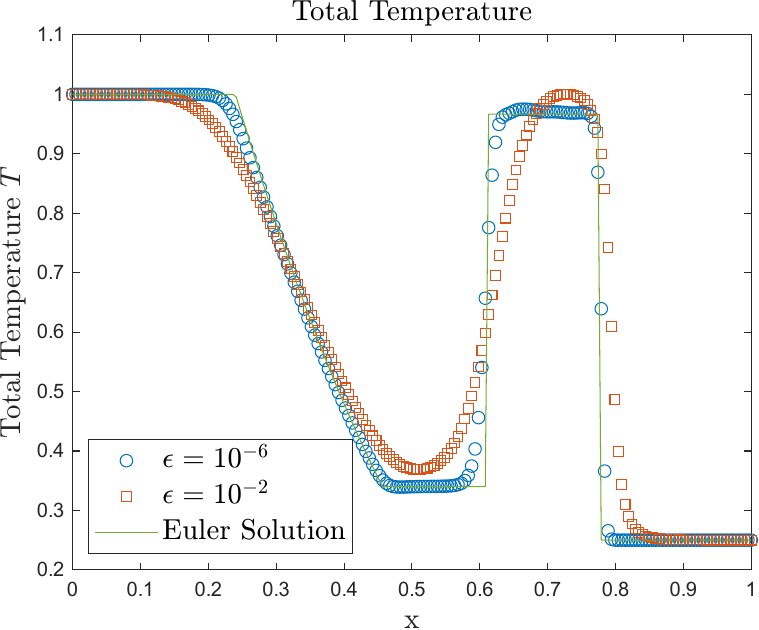}
    }    
    \vspace{0.5cm}
    \subfigure[Velocity]{
        \includegraphics[scale=0.56]{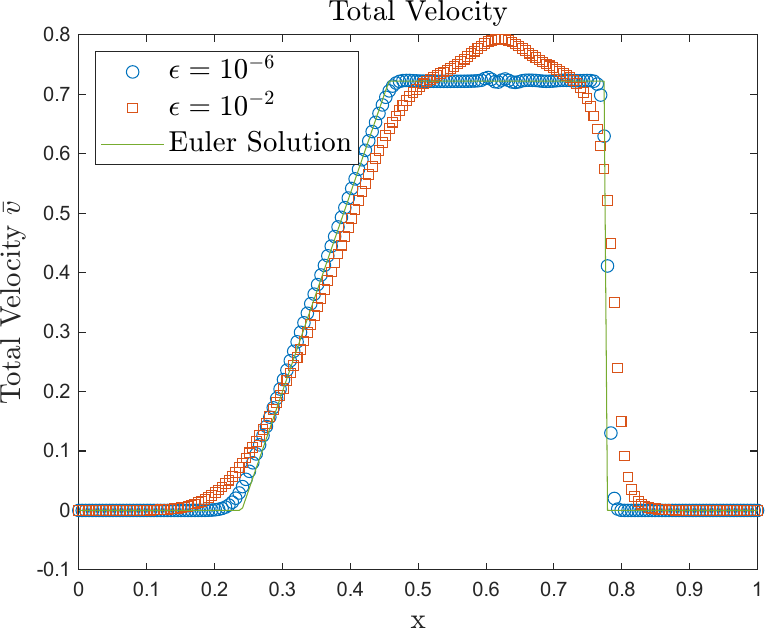}
    }
    \caption{Total density, temperature and velocity for Sod Shock Tube problem (\hyperref[Test1]{\textbf{Test 1}}), for different values of Knudsen number at the final time $t=0.15$ and with unitary mass ratio.} 
    \label{fig_sod_setA}
\end{figure}

\begin{figure}
    \centering
    \subfigure[Density]{
        \includegraphics[scale=0.56]{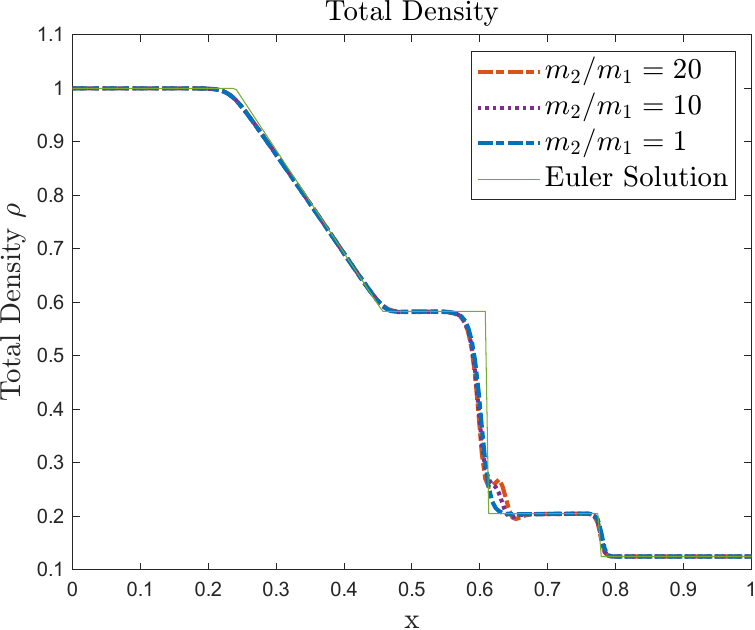}
    }
    \subfigure[Pressure]{
        \includegraphics[scale=0.56]{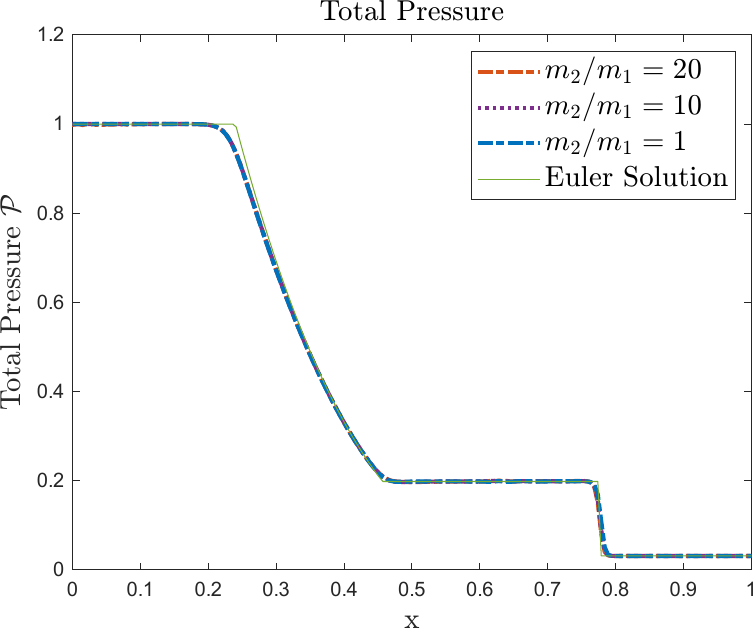}
    }    
    \vspace{0.5cm}
    \subfigure[Velocity]{
        \includegraphics[scale=0.56]{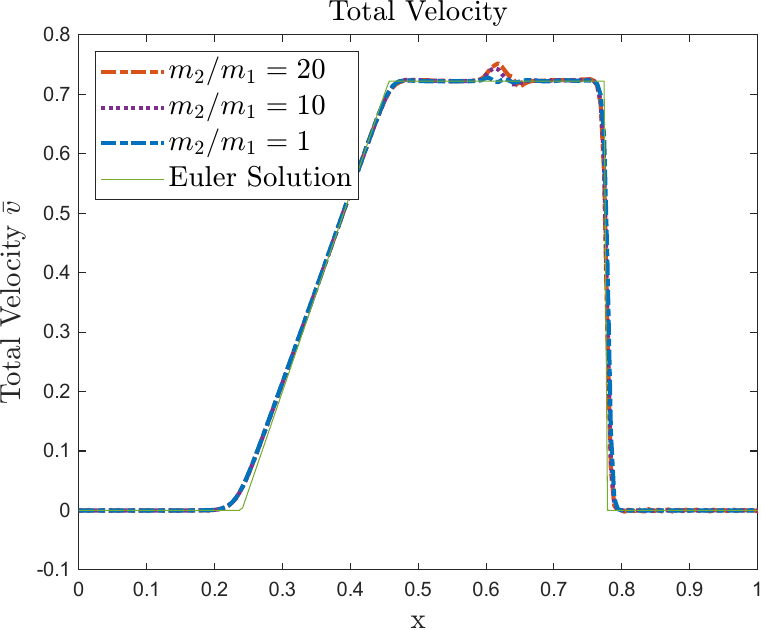}
    }
    \caption{Total Density, Pressure and Velocity in the Sod Shock Tube problem (\hyperref[Test1]{\textbf{Test 1}}), for different values of Mass Ratio ($MR=m_2/m_1$), with a Knudsen number equal to $\varepsilon=10^{-6}$, at the final time $t=0.15$.} 
    \label{fig_sod_setB}
\end{figure}

 A classical benchmark test for the 1D Euler equations is the so-called \emph{Sod shock tube test problem} \cite{tor-2009}: its solution consists of a left rarefaction, a contact and a right shock. Here we have considered its extension to the multi-species configuration using similar initial conditions proposed by \cite{bai-rey-2021}. We consider two species in a 1D domain $[0, L_x]$ where $L_x=1$, discretized using $N_x=200$ cells.\\
In order to initially confine the first gas on the left side of the discontinuity, and the second gas on the right side of the discontinuity we initialise the two distribution functions as Maxwellians whose moments are
\begin{equation}
    \begin{cases}
        \rho_1 = (1-\delta)\rho_L, \  \  \rho_2 = \delta \rho_L, \  \  v_1 = v_2 = v_L, \  \  T_1 = T_2 =\mathcal{P}_L/(\rho_1/m_1), \  \  if\ \  x< 0.5,\\
        \rho_1 = \delta\rho_R, \  \  \rho_2 = (1-\delta) \rho_R, \  \  v_1 = v_2 = v_R, \  \  T_1 = T_2 =\mathcal{P}_R/(\rho_2/m_2), \  \  if\ \  x\geq 0.5,
    \end{cases}
\end{equation}
where $\rho_L=1$, $v_L=0$, $\mathcal{P}_L=1$, $\rho_R=2^{-3}$, $v_R=0$, $\mathcal{P}_R=2^{-5}$.\\
As in \cite{bai-rey-2021}, we set $\delta = 10^{-5}$, because its ideal value ($\delta=0$), would lead to ill-defined macroscopic moments.\\
Each species has its own mass, thus we introduce the mass $m_1$ associated to the first species, and the mass $m_2$ associated to the second species.\\
In all simulations we considered $\Xi = 3$ as specific heat ratio in the computation of the analytical solution.\\
All simulations are executed up to a final integration time equal to $0.15$.\\
We executed two sets of simulations and we compared them with the analytical solution formally obtained when $\varepsilon\to 0$. This reference solution is obtained analytically solving the Riemann problem.
\begin{itemize}
    \item Set A: two simulations with mass ratio equal to 1, and two different values of the Knudsen number, $\varepsilon=10^{-2}$ and $\varepsilon=10^{-6}$. The timestep is equal to $\Delta t = 0.1 dx$. The goal of this test is to check the correct implementation of the integration method and that the scheme is able to detect the correct asymptotic behaviour (Euler equations) when the Knudsen number goes to zero. In both simulations we considered a velocity domain given by $[-L, L]$ with $L=8$ and discretized with 32 points. The density, velocity and temperature of the mixture at the final time are showed in Figure \ref{fig_sod_setA}, from which it is evident that the solution computed with the smallest value of the Knudsen number is converging towards the hydrodynamic (analytic) solution ($\varepsilon\to 0$) and the correctness of the numerical schemes.
    \item Set B: three simulations with different mass ratios: $MR=1$ $MR=10$ and $MR=20$, and with $\varepsilon=10^{-6}$. This is an important test to check that the scheme is properly able to numerically integrate those configurations where the two species have different masses. The velocity space of the simulation with $MR=1$ is given by $[-L, L]$ with $L=8$ and it is discretized with 32 points, and the timestep is equal to $\Delta t = 0.1 dx$; while for the simulations with $MR=10$ and $MR=20$ we considered $L=40$ and 200 points, and the timestep is equal to $\Delta t = 0.025 dx$. The densities, velocities and pressure of the mixture are showed in Figure \ref{fig_sod_setB}. It is evident that the code is able to correctly capture the correct asymptotic behaviour, even with high mass ratios. The bumps that can be seen in the numerical solutions have been already detected and discussed by other authors, like \cite{bai-rey-2021}: they are due to the effect of high mass ratio in proximity of the contact discontinuity between the two species. The choice of the velocity box can affect the magnitude of these bumps \cite{bai-rey-2021}.
\end{itemize}

\subsubsection{Test 2: Kelvin-Helmholtz instability}
\label{Test2}
\begin{figure}
    \centering
    \subfigure[Species 1, $\rho_1$]{
        \includegraphics[scale=0.49]{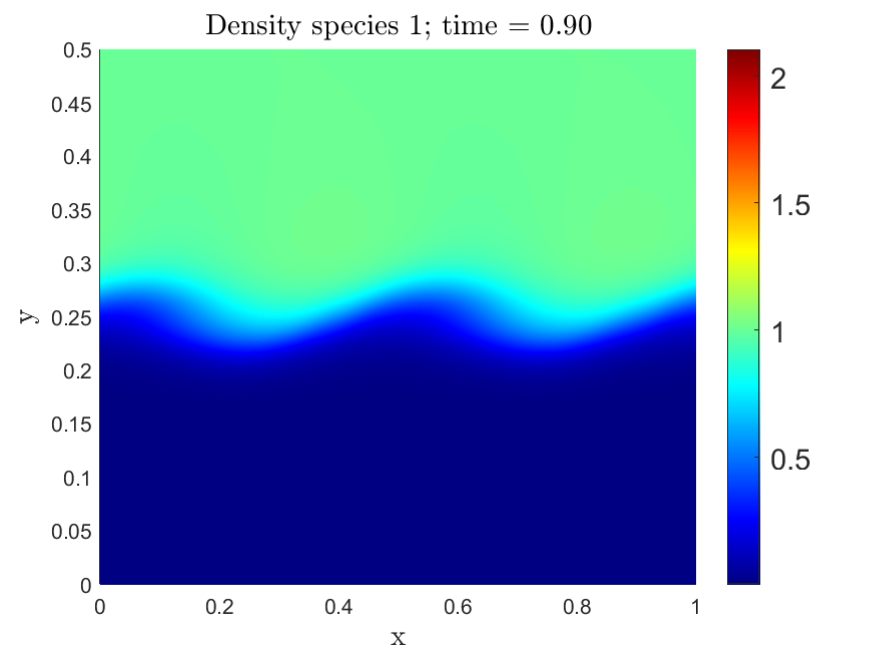}
    }
    \subfigure[Species 2, $\rho_2$]{
        \includegraphics[scale=0.49]{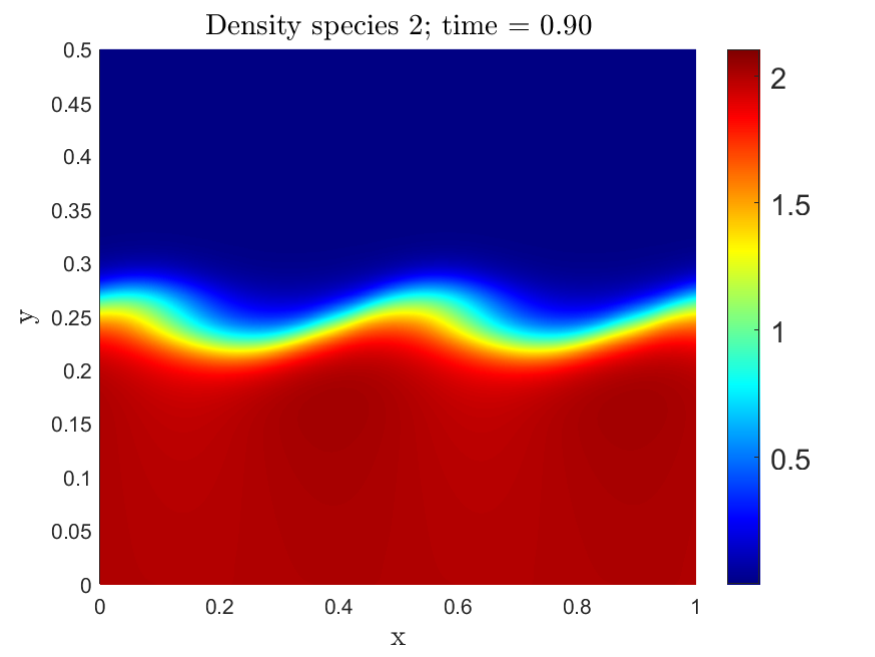}
    }
    \hfill
    \centering
    \subfigure[Species 1, $\rho_1$]{
        \includegraphics[scale=0.49]{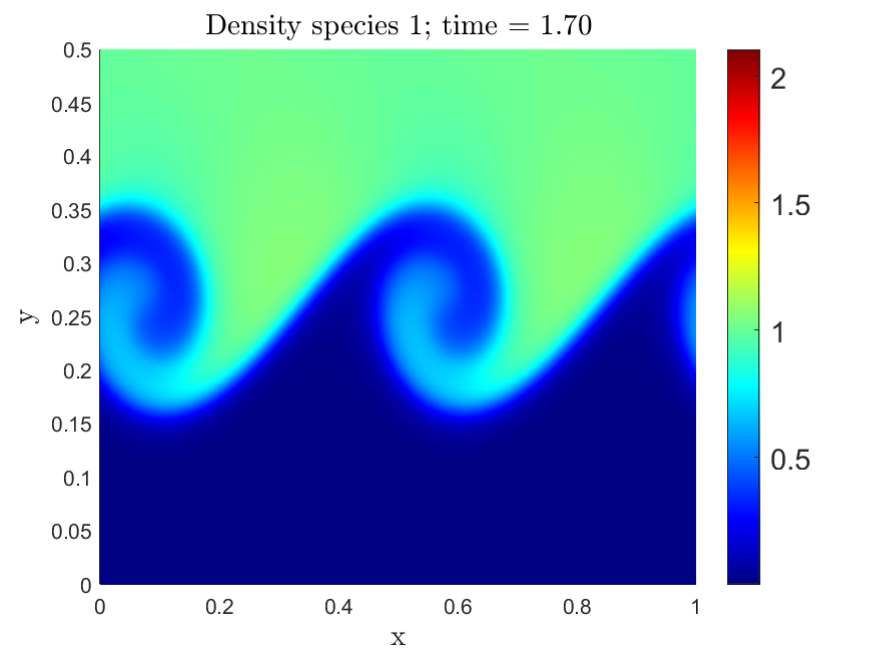}
    }
    \subfigure[Species 2, $\rho_2$]{
        \includegraphics[scale=0.49]{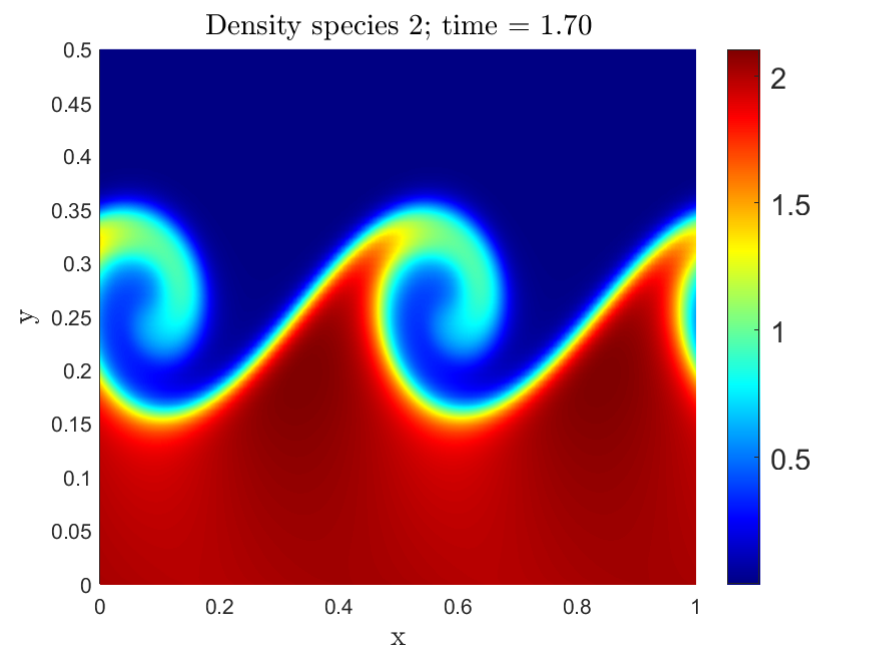}
    }  
    
    \hfill
    \centering
    \subfigure[Species 1, $\rho_1$]{
        \includegraphics[scale=0.49]{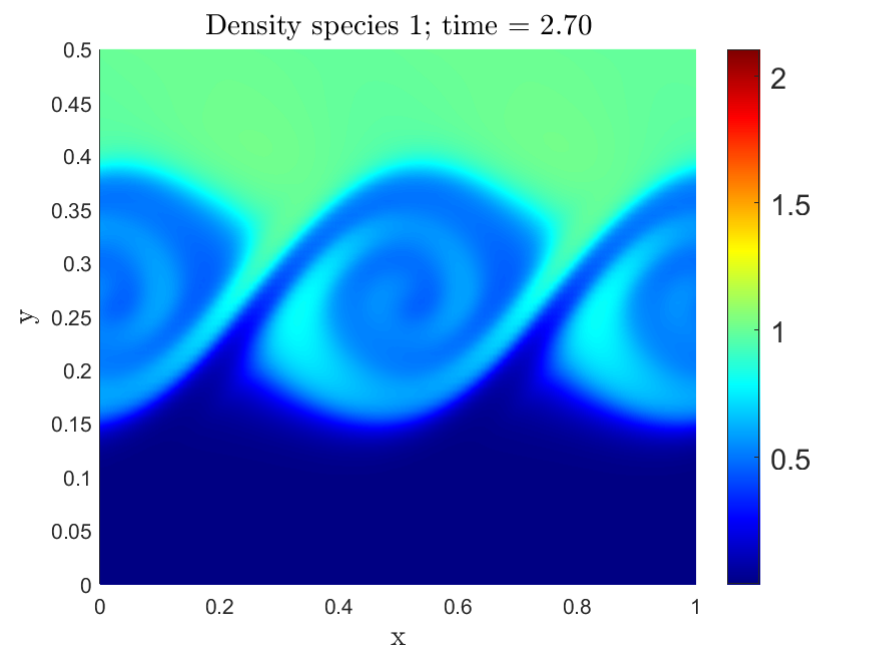}
    }
    \subfigure[Species 2, $\rho_2$]{
        \includegraphics[scale=0.49]{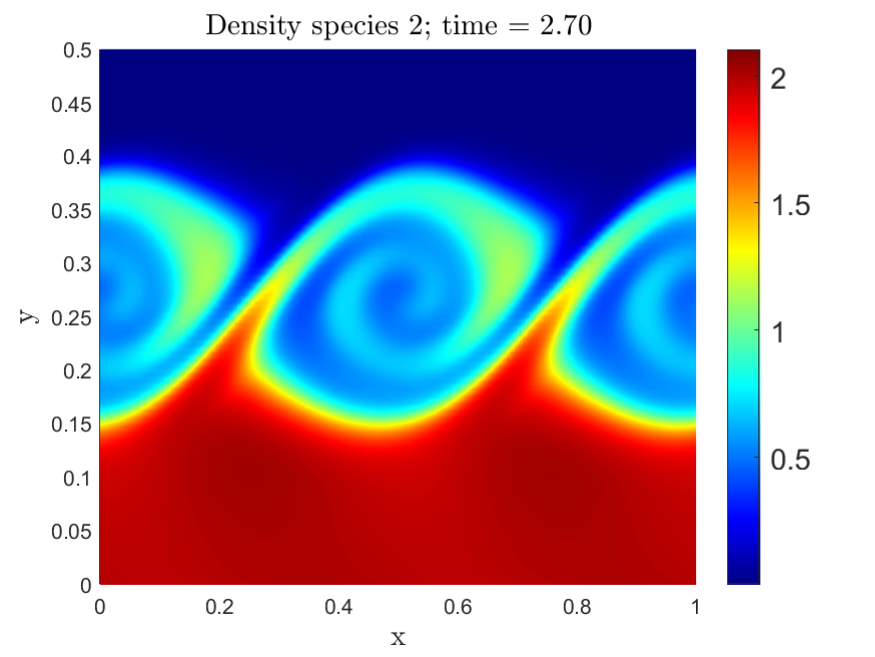}
    }
    \caption{Densities of the two species at $t=0.9$, $t=1.7$ and $t=2.7$ for the Kelvin-Helmholtz instability (\hyperref[Test2]{\textbf{Test 2}}).
    } 
    \label{fig_KH_fig1}
\end{figure}

The Kelvin-Helmholtz instability is a phenomenon where vortices arise at the interface of two fluids moving at different horizontal speeds and with a perturbative vertical sinusoidally spatially distributed velocity.\\
The initial conditions are inspired from \cite{mcn-lyr-pas-2012, mel-rey-sam-2017} and have been adapted for Gas-Mixtures by Bailo and Rey in \cite{bai-rey-2021}. This is an important test to check the implementations of our 2D integration schemes for our kinetic model.\\
We consider two gases with mass ratio $m_2/m_1=2$ in a domain $[0, 1]\times [0, 0.5]$, periodic along x-direction, and discretized using $N_x=150$ and $N_y=75$ cells. The velocity domain $[-8, 8]^2$ is discretized using 32 points in both directions. The Knudsen number is equal to $\varepsilon=5\times10^{-5}$. The two gases are initially separated along the interface surface placed at $y=L_y/2$ and they are initialized as Maxwellians whose moments correspond to the following values:
\begin{align}
    y> L_y/2: 
    \begin{cases}
        \rho_1=\tilde{\rho}_1(1-\delta), \  \  v_1=(0.5, 0.01\sin{(4\pi x)}), \  \  T_1=(\tilde{\rho}_1/m_1)^{-1},\\
        \rho_2=\tilde{\rho}_2\delta , \  \  v_2=(0.5, 0.01\sin{(4\pi x)}), \  \  T_2=T_1=(\tilde{\rho}_1/m_1)^{-1},\\
    \end{cases}\\
    y\leq L_y/2: 
    \begin{cases}
        \rho_1=\tilde{\rho}_1\delta, \  \  v_1=(-0.5, 0.01\sin{(4\pi x)}), \  \  T_1=T_2=(\tilde{\rho}_2/m_2)^{-1},\\
        \rho_2=\tilde{\rho}_2(1-\delta) , \  \  v_2=(-0.5, 0.01\sin{(4\pi x)}), \  \  T_2=(\tilde{\rho}_2/m_2)^{-1},\\
    \end{cases}
\end{align}
where $\delta=10^{-5}$, $\tilde{\rho}_1=1$ and $\tilde{\rho}_2=2$.\\
The densities of the two species for three different times $t=0.9$, $t=1.7$ and $t=2.70$ are plotted in Figure \ref{fig_KH_fig1}.
At the final timestep, it is clear the formation of the vortices, consistently with \cite{bai-rey-2021}.

\subsubsection{Test 3: Flow around a cylinder}
\label{Test3}
\begin{figure}
    \centering
    \subfigure[Species 1, $\rho_1$]{
        \includegraphics[scale=0.49]{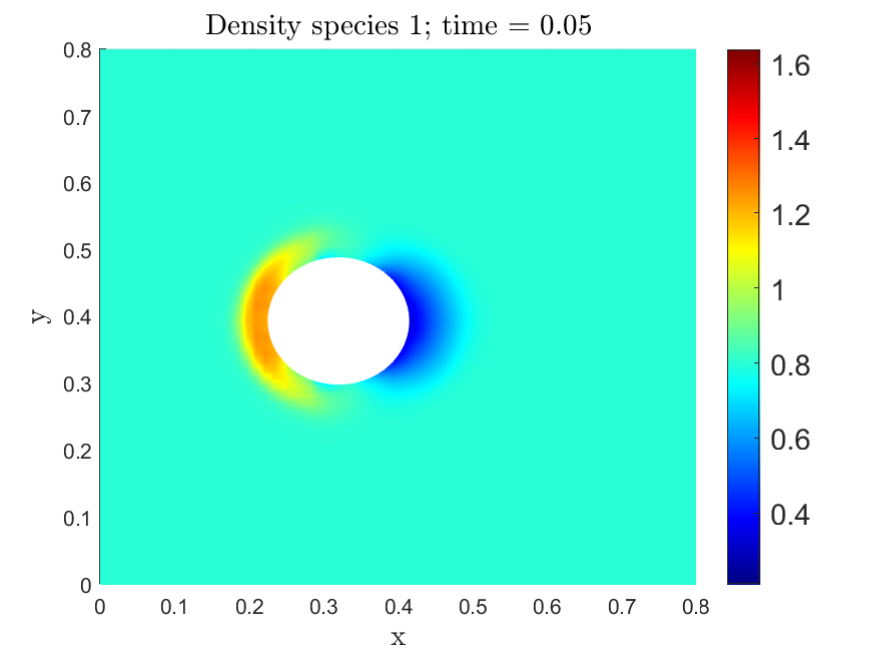}
    }
    \subfigure[Species 2, $\rho_2$]{
        \includegraphics[scale=0.49]{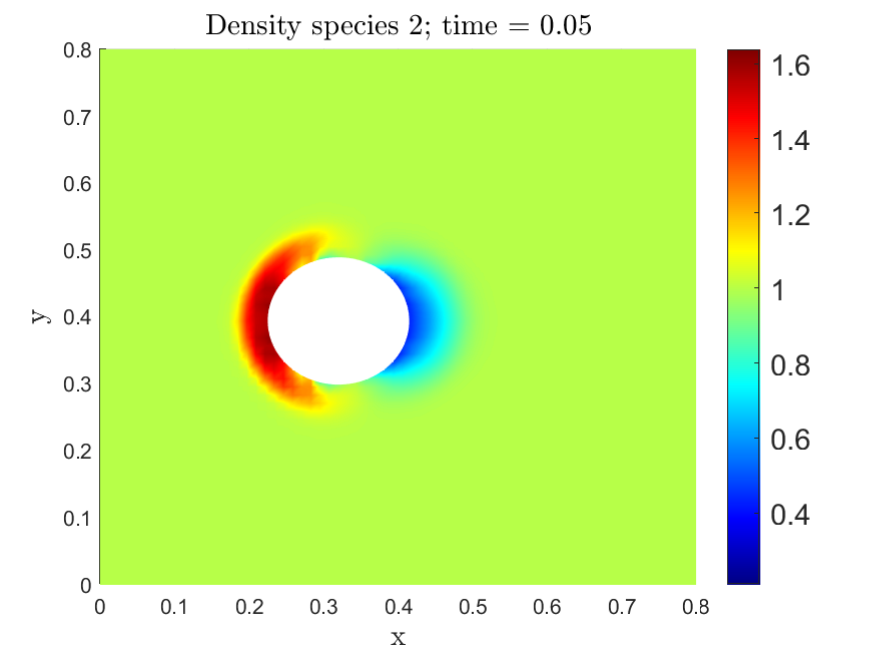}
    }
    \vspace{0.5cm}
    \subfigure[Species 1, $\rho_1$]{
        \includegraphics[scale=0.49]{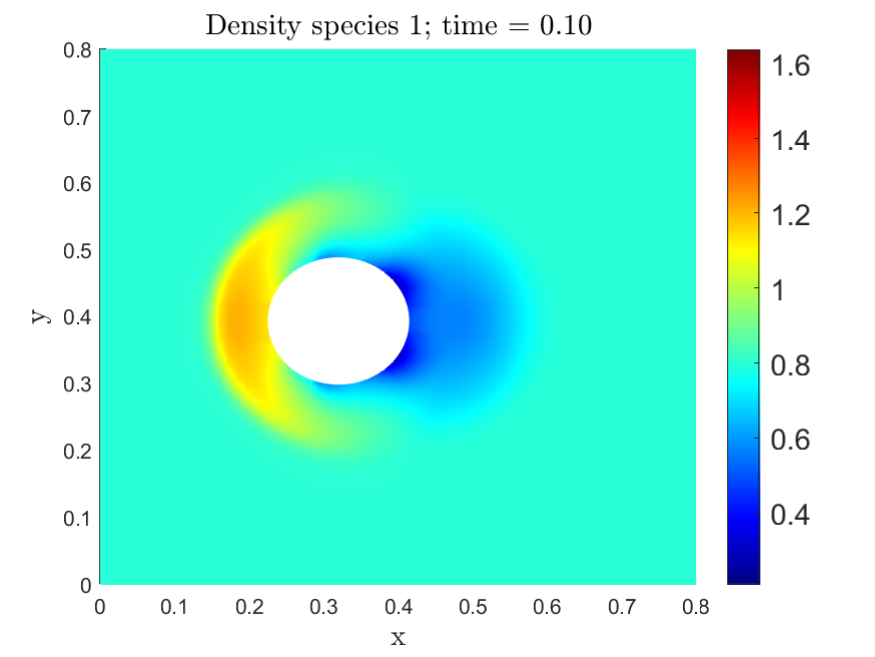}
    }
    \subfigure[Species 2, $\rho_2$]{
        \includegraphics[scale=0.49]{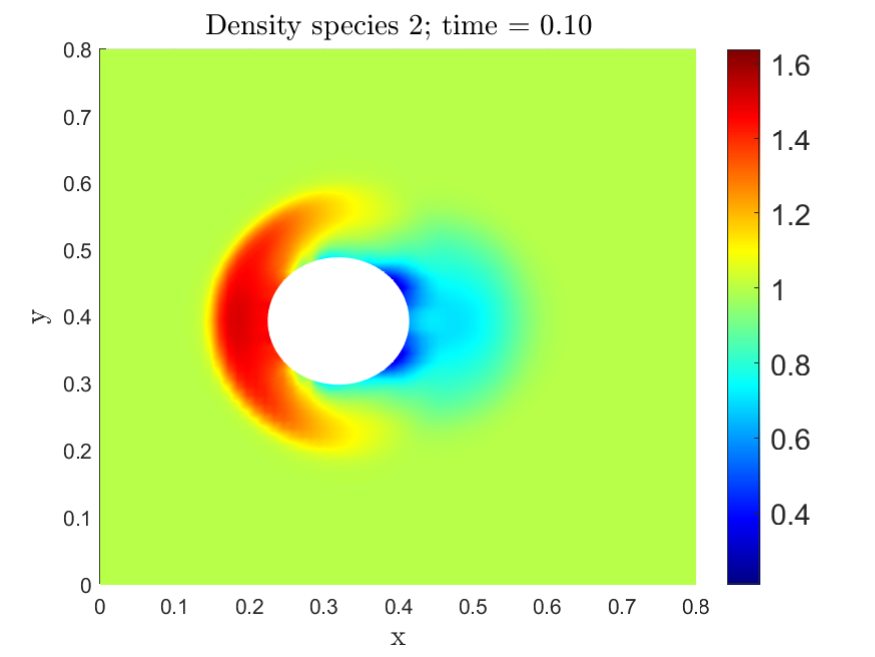}
    }
    \vspace{0.5cm}
    \subfigure[Species 1, $\rho_1$]{
        \includegraphics[scale=0.49]{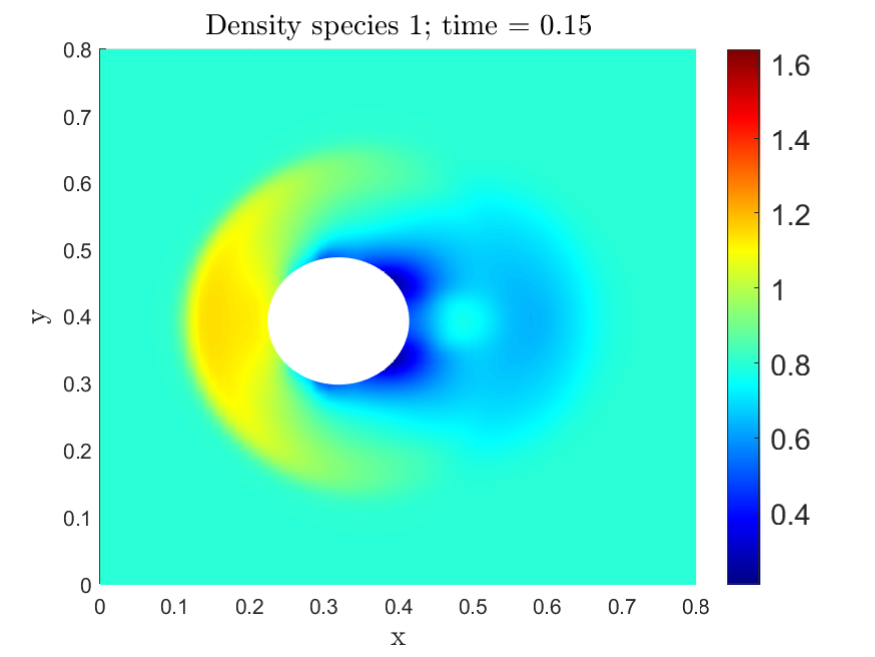}
    }
    \subfigure[Species 2, $\rho_2$]{
        \includegraphics[scale=0.49]{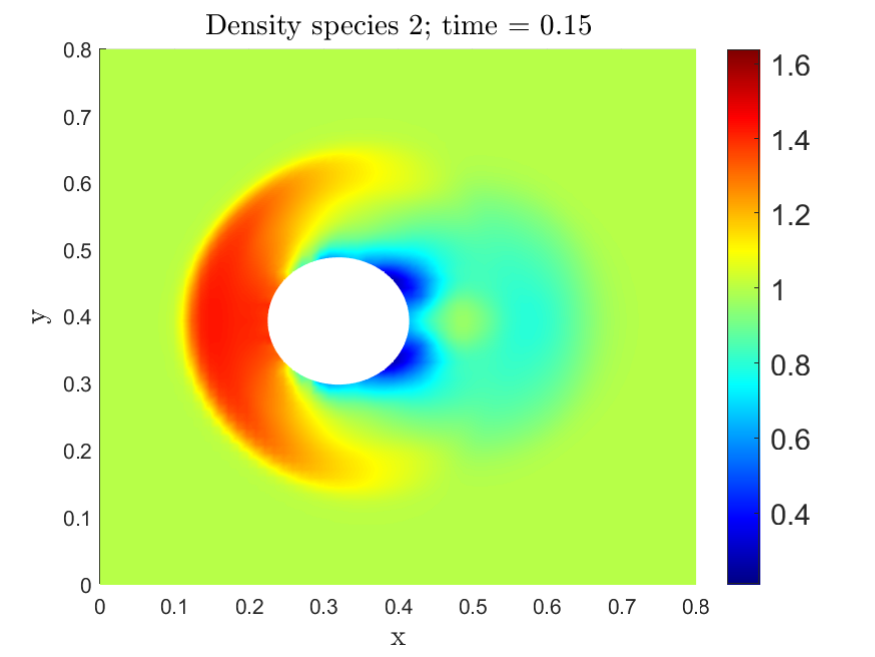}
    }
    \caption{Densities evolution of the two species at $t=0.05$, $t=0.10$ and $t=0.15$ for the flow around a circular shape (\hyperref[Test3]{\textbf{Test 3}}).} 
    \label{fig_cylinder_fig1}
\end{figure}

The aim of this test is to verify that the numerical scheme is able to handle complex geometries, like a flow over a round object, similarly to \cite{fil-jin-2011, yan-hua-1995}. We consider a 2D domain of size $[0, L_x]\times[0, L_y]$ where $L_x=L_y=0.8$, and discretized with $66$ cells in both directions. The velocity space is a domain equal to $[-L,L]^2$ where $L=8$ and it is discretized using $32$ points in both directions. The Knudsen number is equal to $10^{-6}$ and we consider a mass ratio $MR=m_2/m_1=5$ ($m_1=1$ and $m_2=5$). The round object is modelled considering a round shape centred in the cell of indices $[4L_x/10, N_y/2]$ and whose radius is equal to $R=L_x/10$ in physical coordinates, such that the integrator is never executed inside this round shape. We implement specular boundary conditions \cite{cer-ill-pul-1994} on the walls of the object. The distribution functions for both species are initialized like a Maxwellian whose moments outside of the round object are given by
\begin{equation}
    \begin{cases}
        \rho_1 = 0.8, \  \  v_1 = (0.8, 0), \  \  T_1 = 1,\\
        \rho_2 = 1.0, \  \  v_2 = (0.8, 0), \  \  T_2 = 1\\
    \end{cases}.
\end{equation}
In Figure \ref{fig_cylinder_fig1} they are plotted the densities of the two species for three different times $t=0.05$ and $t=0.10$ and $t=0.15$, from which it is evident that the code is able to handle complex geometries. The timestep is approximately equal to $\Delta t=7.7\times 10^{-4}$. 

From this test it is also clear that the densities of the two species are different one from each other, due to their different masses and different initial values: thus it is remarkable the importance of using a multi-species framework rather than the single species one.

\subsubsection{Test 4: Comparison between BGK, ES-BGK and full Boltzmann operator}
\label{Test4}
\begin{figure}
    \centering
    \subfigure[$\varepsilon=10^{-1}$]{
        \includegraphics[scale=0.55]{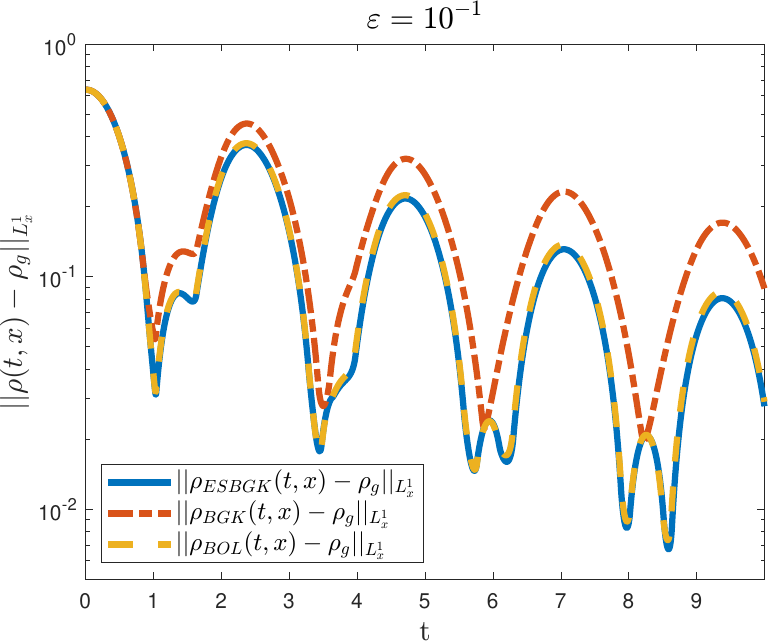}
        \label{fig:enter-test5_1} 
    }
    \subfigure[$\varepsilon=10^{-4}$]{
        \includegraphics[scale=0.55]{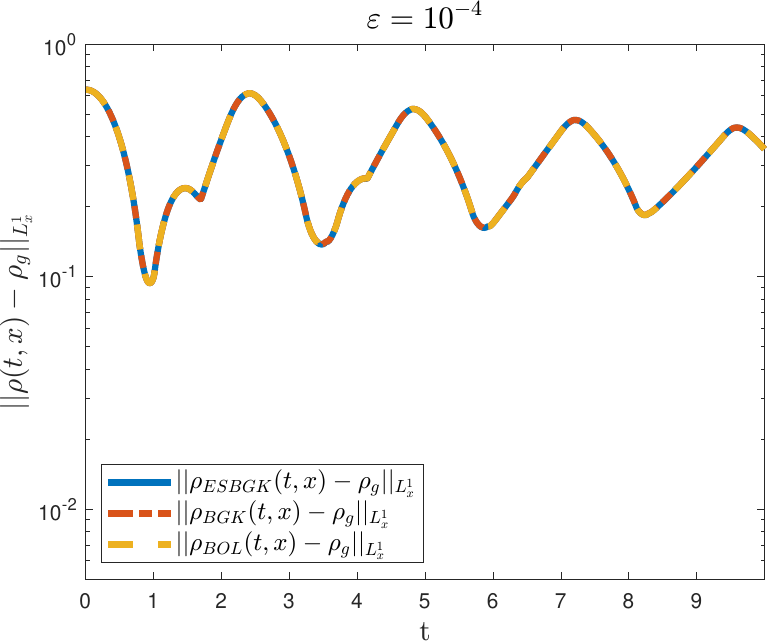}
        \label{fig:enter-test5_2} 
    }
    \caption{Evolution of the quantity $||\rho(t, x)-\rho_g||_{L^1_x}^{\text{norm}}$ (\hyperref[Test4]{\textbf{Test 4}}) obtained using three different operators: multi-species ES-BGK, multi-species BGK and single species spectral full Boltzmann, for two different values of Knudsen number $\varepsilon$. The multi-species solvers are used to evolve the distribution functions $f_1$ and $f_2$ of two species that have the same mass and the same cross section, while the single-species Boltzmann operator is used to evolve the single distribution $f$ initialized as the sum of the two species $f(t=0)=f_1(t=0)+f_2(t=0)$. Then thanks to the indifferentiability principle it is possible to compare the evolution of the sum of the two distribution functions $f_1(t)+f_2(t)$ with the the evolution of $f(t)$ for every time $t$.} 
    \label{fig_comparison}
\end{figure}

The aim of this test, inspired by \cite{fil-jin-2011}, is to show the ability of our multi-species ES-BGK scheme to reproduce the correct behaviour of the full Boltzmann operator.\\
For simplicity we consider the case of two species with same unitary masses $m_1=m_2=m=1$. In what follows we strongly use this assumption, that it is also really useful, because allows us to check the indifferentiability principle \cite{bru-2015} of our scheme: when all the species have the same masses and their crosses section are equal, then the system of equations reduces to a single one by adding the distribution functions.\\
Here we will follow and summarize \cite{fil-jin-2011} in presenting this test performed in a 1D space domain $[-L_x, L_x]$, with $L_x=1$, and 2D velocity space. We will consider two species of same masses, and the distribution functions initialized as $f_1(t=0)=f_2(t=0)$ and evolved using the scheme presented in this paper.\\
Considering a smooth solution with periodic boundary conditions, in a 1D space domain $[-1, 1]$, from the conservation laws for mass, momentum and energy, one can determine the stationary state of the distribution function $f=f_1+f_2$, which is the normalized global Maxwellian distribution \cite{fil-jin-2011}
\begin{equation}
    \mathcal{M}_g(v) = \frac{\rho_g}{2\pi T_g}\exp{\left(-\frac{|v|^2}{2T_g}\right)},
\end{equation}
where $\rho_g$ and $T_g$ to be determined, under the assumption that
\begin{equation}
    \frac{1}{2}\int_{-1}^{1}\int_{\RR^2}f_1(t=0) \  v \  dv \  dx = \frac{1}{2} \int_{-1}^{1}\int_{\RR^2}f_2(t=0) \  v \  dv \  dx = 0.
\end{equation}
Under some boundness assumption, it is expected the distribution function $f$ converges to the global Maxwellian $\mathcal{M}_g$ as $t\to \infty$. One expects that the if the solution is smooth enough then \cite{fil-jin-2011}
\begin{equation}
    ||f(t)-\mathcal{M}_g||=\mathcal{O}(t^{-\infty}),
\end{equation}
meaning that the solution converges almost exponentially fast to the global equilibrium.\\
In order to check this convergence behaviour we study the evolution of the following quantity
\begin{equation}
    ||\rho(t, x)-\rho_g||_{L^1_x}=
    \int_{-1}^{1}|\rho(t, x)-\rho_g| \ dx,
\end{equation}
where $\rho$ is the total mass density of the distribution function $f$ (where $f=f_1+f_2$) and
\begin{equation}
    \rho_g = \frac{1}{2}\int_{-1}^{1}\int_{\RR^2}\left(f_1(t=0, x, v) + f_2(t=0, x, v) \right) \ dv \ dx.
\end{equation}
In particular we will compare the evolution of this quantity obtained from three different approaches:
\begin{itemize}
    \item multi-species ES-BGK model presented in this paper by setting $\nu=-1$ in our scheme,
    \item multi-species BGK model obtained by putting $\nu=0$ in our scheme,
    \item spectral full Boltzmann operator for single species \cite{fil-mou-par-2006}.
\end{itemize}
We use the first two integrators to evolve the multi-species system given by $f_1$ and $f_2$, initialized using the same distribution function and having the same unitary mass; and we use the third one to evolve $f$ which is initialized as $f(t=0)=f_1(t=0)+f_2(t=0)$. Indeed thanks to the indifferentiability principle the total density $\rho$ obtained from the multi-species model ($\rho=\rho_1+\rho_2$) and the total density $\rho$ obtained from the single-species such defined, are the same.\\
In order to obtain the correct full Boltzmann Prandtl number with the ES-BGK operator we set $\nu=-1$ \cite{fil-jin-2011}. By equating the viscosity of the Boltzmann operator with the one of the ES-BGK operator (see equation \eqref{eqViscosity}) one gets a constraint on the parameter $\lambda$ \cite{fil-jin-2011}, that is $\lambda=\frac{\mathcal{P}}{1-\nu} \frac{A_2 3 \pi}{\sqrt{2}T}$, where $A_2=0.436$. We use this equation to set the value of $\lambda$ for the ES-BGK operator (by setting $\nu=-1$) and for the BGK operator (by setting $\nu=0$).
We initialize the two distribution functions as
\begin{equation}
    f_1(t=0, x,v)=f_2(t=0, x, v)=\frac{1+A_0 \sin{(\pi x)}}{4\left(2\pi T_0\right)} \left[\exp{\left(-\frac{|v-u_0|^2}{2T_0}\right)}+\exp{\left(-\frac{|v+u_0|^2}{2T_0}\right)}\right],
\end{equation}
where $A_0=0.5$, $T_0=0.125$ and $u_0=(1/2, 1/2)$. We update $f_1$ and $f_2$ using the multi-species ES-BGK operator, and the multi-species BGK operator, while we update $f$ (initialized as $f(t=0)=f_1(t=0)+f_2(t=0)$) using the full Boltzmann operator. We consider $L_x=1$  with a periodic physical space $[-L_x, L_x]$ discretized using 100 cells, and a velocity box $[-8,8]^2$ discretized using 64 points in both directions.\\
In all simulations the space reconstruction is performed using the already discussed CWENO3 scheme. We perform the test using two different values of the Knudsen number $\varepsilon=10^{-1}$ and $\varepsilon=10^{-4}$. For the sake of simplicity, the time integration of the ES-BGK and BGK operators is performed using the first order AP IMEX scheme discussed in this paper, while for the full Boltzmann operator we use the first order IMEX AP scheme proposed in \cite{fil-jin-2010}. The timestep is equal to $dt=0.1 dx$ in all simulations, and the solutions are evolved up to a final time $t_F=10$.\\
The curves $||\rho(t, x)_{ESBGK}-\rho_g||_{L^1_x}$, $||\rho(t, x)_{BGK}-\rho_g||_{L^1_x}$, $||\rho(t, x)_{BOL}-\rho_g||_{L^1_x}$, reported in Figure \ref{fig_comparison}, are obtained respectively using the multi-species ES-BGK operator, the multi-species BGK operator and the single species spectral full Boltzmann operator. An important indicator of the correctness our results, are the curves in Figure \ref{fig_comparison} obtained with $\varepsilon=0.1$, which are in perfect agreement with those presented in \cite{fil-jin-2011}. It is evident that the ES-BGK operator is able to correctly match and capture the behaviour of the curve obtained using the full Boltzmann operator even in the rarefied regime $\varepsilon=0.1$, while the BGK operator only in the collision dominated regime $\varepsilon=10^{-4}$. This test clearly shows the importance of using the ES-BGK operator rather than the BGK operator in the rarefied regime.\\

\subsubsection{Discussion on the IMEX scheme computational performance}

We now briefly discuss the speed-up offered by the proposed IMEX scheme when compared to a standard fully explicit integrator, using our numerical tests as a representative example. More details can be found in \cite{bos-par-rus-2024}.
The first clear advantage of the IMEX scheme is the already discussed Asymptotic-Preserving property, which ensures that, for sufficiently small Knudsen numbers, the method becomes an exact solver for the multi-species Euler and Navier–Stokes limits (see Section~\ref{APsection}).
Moreover, the IMEX structure guarantees that the method can be implemented at the cost of a fully explicit scheme but with a timestep  constrained only by the CFL condition associated with the transport term, and not by the value of the Knudsen number. In contrast, in a standard fully explicit solver the timestep must satisfy
\[
\Delta t^{\text{Explicit}} \lesssim \min\{\varepsilon,\ \Delta x / v_{\max} \},
\]
where $\varepsilon$ is the Knudsen number, $v_{\max}$ the maximal velocity in the velocity domain, and $\Delta x$ the spatial grid size.

To make this concrete, consider \hyperref[Test3]{\textbf{Test 3}} (flow around a cylinder).  
In this simulation we used:
\[
\varepsilon = 10^{-6},\qquad \Delta x = 0.0123,\qquad v_{\max} = 8.
\]
Hence,
\[
\min\{\varepsilon, \Delta x / v_{\max}\} = \varepsilon = 10^{-6},
\]
so a fully explicit solver would be forced to use a timestep of order $\Delta t^{\text{Explicit}} \approx 10^{-6}$.

In contrast, our IMEX scheme is constrained only by $\Delta x / v_{\max}$, yielding a timestep of approximately
\[
\Delta t^{\text{IMEX}} = 7.7\times 10^{-4}.
\]
Therefore, the effective speed-up is
\[
\frac{\Delta t^{\text{IMEX}}}{\Delta t^{\text{Explicit}}} \approx 770,
\]
which is extremely significant.

Of course, the precise value of the speed-up depends on both the Knudsen number and the CFL constraint, which are problem-dependent. Remarkably, in the worst-case scenario (when $\varepsilon > \Delta x / v_{\max}$), the maximal allowable IMEX timestep coincides with that of the explicit scheme. Thus, IMEX schemes are \emph{always} at least as good as explicit ones, and often dramatically better.

This qualitative discussion highlights the importance of developing IMEX methods for accelerating stiff kinetic simulations.
Finally, let us mention that further reductions in computational time could be achieved by introducing a local velocity grid adaptation strategy, such as the one proposed in \cite{ber-iol-pup-2014}.

\section{Conclusion}\label{conclusions}
We have developed a high-order Asymptotic-Preserving IMEX scheme for the ES-BGK model of gas mixtures. The construction extends IMEX strategies previously available for single-species BGK models to the multi-species ES-BGK collision operator and is formulated in such a way that the time discretization is explicitly implementable and uniformly stable with respect to the Knudsen number. In the fluid regime, the scheme reduces to a consistent and high-order accurate discretization of the limiting macroscopic equations for the mixture, while in the rarefied regime it remains a genuine kinetic solver.

The proposed method has been embedded in a high-order numerical framework that combines third-order IMEX Runge--Kutta time integration, CWENO3 spatial reconstruction, conservation of macroscopic moments in the discrete velocity space, and a multithreaded implementation. The numerical tests confirm the expected accuracy of the scheme in both kinetic and fluid regimes, as well as a smooth transition between them, and illustrate its capability to handle an arbitrary number of species without loss of stability.

Several directions for future work are natural. On the modeling side, one can consider the extension of the present framework to more complex collisional models, and to mixtures in plasma physics. On the numerical side, it would be of interest to couple the present IMEX strategy with different high-order spatial discretization and unstructured meshes, and to investigate its performance in multi-scale applications where kinetic and fluid regions coexist in a more complex geometry.

\section*{Acknowledgments}
The authors received funding from the European Union's Horizon Europe research and innovation program under the Marie Skłodowska-Curie Doctoral Network DataHyking (Grant No. 101072546). 
DC and TR are also supported by the French government, through the UniCA$_{JEDI}$ Investments in the Future project managed by the National Research Agency (ANR) with the reference number ANR-15-IDEX-01.
DC acknowledges the support of the Mathematics Department of the University of Ferrara through access to their computational resources. The research of LP has been supported by the Royal Society under the Wolfson Fellowship “Uncertainty quantification, data-driven simulations and learning of multiscale complex systems governed by PDEs”. LP also acknowledges the  partial support by Fondo Italiano per la Scienza (FIS2023-01334) advanced grant ADAMUS. This work has been written within the activities of GNCS group of INdAM (Italian National Institute of High Mathematics).  
\printbibliography
\end{document}